\documentclass{amsart}
\usepackage{latexsym}
\usepackage{rotating}
\usepackage{epsfig}

\font\got=cmss10

\newtheorem{Thm}{Theorem}[section]
\newtheorem{Lm}[Thm]{Lemma}

\newtheorem{lemma}[Thm]{Lemma}
\newtheorem{Prop}[Thm]{Proposition}
\newtheorem{Cor}[Thm]{Corollary}

\theoremstyle{definition}
\newtheorem{defn}[Thm]{Definition}

\theoremstyle{definition}
\newtheorem{Remark}[Thm]{Remark}
\newtheorem{Example}[Thm]{Example}

\renewenvironment{proof}{\noindent {\em Proof. }}{\hfill $\Box$}

\DeclareMathOperator{\ad}{ad}
\def\alp{\alpha}
\DeclareMathOperator{\Aut}{Aut}
\def\gam{\gamma}

\renewcommand{\phi}{\varphi}
\def\spe{extremal}
\DeclareMathOperator{\tr}{tr}

\newcommand{\erz} [1] {\mbox{$\langle #1 \rangle$}}
\newcommand{\ov}{\overline}

\newcommand{\sli}{{\mathfrak{sl}}}

\newcommand{\spa}{\operatorname{span}}

\newcommand{\Rad}{\operatorname{Rad}}
\newcommand{\NilRad}{\operatorname{NilRad}}
\newcommand{\SanRad}{\operatorname{SanRad}}
\newcommand{\ch}{\operatorname{char}}
\newcommand{\wt}{\widetilde}
\newcommand{\eps}{\varepsilon}
\newcommand{\im}{{\rm Im}\,}

\newcommand{\ek}{{\mathfrak{e}}}
\newcommand{\fk}{{\mathfrak{f}}}
\newcommand{\g}{{\mathfrak{g}}}
\newcommand{\hk}{{\mathfrak{h}}}
\newcommand{\nk}{{\mathfrak{n}}}
\newcommand{\Bc}{{\mathcal{B}}}
\newcommand{\CC}{{\mathbb{C}}}

\newcommand{\NN}{{\mathbb{N}}}
\newcommand{\ZZ}{{\mathbb{Z}}}
\newcommand{\Ec}{{\mathcal{E}}}
\newcommand{\Dc}{{\mathcal{D}}}
\newcommand{\Lc}{{\mathcal{L}}}
\newcommand{\Rc}{{\mathcal{R}}}

\begin{document}
\title{Lie algebras generated by extremal elements}
\author{Arjeh M.~Cohen, Anja Steinbach, Rosane Ushirobira, David Wales}
\address{}
\date{\today}

\begin{abstract}
We study Lie algebras generated by extremal elements (i.e., elements
spanning inner ideals of $L$) over a field of characteristic distinct
from 2. We prove that any Lie algebra generated by a finite number of
extremal elements is finite dimensional. The minimal number of
extremal generators for the Lie algebras of type $A_n$
$(n\ge1)$, $B_n$ $(n\ge3)$, $C_n$ $(n\ge2)$, $D_n$ $(n\ge4)$, $E_n$
$(n=6,7,8)$, $F_4$ and $G_2$ are shown to be $n+1$, $n+1$, $2n$, $n$,
5, 5, and 4 in the respective cases. These results are related to group
theoretic ones for the corresponding Chevalley groups.
\end{abstract}

\maketitle

\section{Introduction}

Let $k$ be a field of characteristic not $2$, and $L$ a Lie algebra
over $k$. We study the r\^ole of {\sl extremal\/} elements in $L$,
that is, those $x\in L$ with $[x,[x,L]]\subseteq kx$.  Since an inner
ideal of $L$ is by definition (\cite{Ben}) 
a linear subspace $I$ of $L$ such that
$[I,[I,L]]\subseteq I$, this amounts to $kx$ being an inner ideal.  By
$\Ec$ or, if necessary to express dependence on $L$, by $\Ec(L)$, we
shall denote the set of all nonzero extremal elements of $L$. 

We are mostly interested here in Lie algebras generated by extremal
elements.  The motivation stems from the fact that long root elements
are extremal in Lie algebras of Chevalley type (i.e., those Lie
algebras over $k$ that are given by the multiplication table of a
Chevalley basis, coming from a simple Lie algebra in characteristic
0).  They were used by Chernousov \cite{Che} in his proof of the Hasse
principle for $E_8$. The associated root groups were studied in a more
abstract group theoretic setting in \cite{Tim2}. Sandwiches, that is,
elements $x\in L$ with $[x,[x,L]] = 0$ are extremal elements of a
special kind; they are prominent in the classification of finite
dimensional simple modular Lie algebras over algebraically closed
fields of characteristics 5 and 7 and could well be useful in a similar way for
other positive characteristics, cf.\ \cite{PrStr}.

We shall be particularly concerned with Lie algebras which are
generated by a finite set of extremal elements.  In \S\ref{finsec}, we
show that in this case the Lie algebra has finite dimension.  For this
we use the work of Zel'manov and Kostrikin \cite {ZK}. It follows from their
results that there is a universal Lie algebra $\Lc_r$ generated by a
finite number of sandwich elements $x_1,\ldots,x_r\in \Lc_r$; it is
nilpotent and of finite dimension.  This leads to the question of a
description of these universal Lie algebras.  Some information about
these for small values of $r$ is given.

Using functions arising from the definition of extremal elements, we are able to
define an associative bilinear form $f$ on $L$.  In particular, if
$x\in\Ec$ and $y\in L$, then $[x,[x,y]]$ is a multiple of $x$ and so
we can set $[x,[x,y]]=f(x,y)x$. In \S\ref{bilsec}, we discuss the
connection of $f$ with the Killing form, and give some information
about the radical of this form.  In \S\ref{mingensec}, we determine
the minimal number of generators by extremal elements for the 
Lie algebras of Chevalley type.  Exponentiation by $\ad_x$ for $x$ an extremal
element gives an automorphism of $L$.  The group generated
by these automorphisms satisfies some conditions defined by Timmesfeld
\cite{Tim2}, which we discuss in \S\ref{rootgrpsec}.

We gratefully acknowledge some extremely useful discussions with Gabor Ivanyos
as well as computer verifications of our results on $\Lc_5$
by Willem de Graaf.

\section{General properties of \spe{} elements}
\label{genprop}
\label{3}
Throughout the remainder of this paper,
$k$ is a field of characteristic distinct from 2 and $L$ is a Lie algebra over $k$.
By linearity of $[x,[x,y]]$ in $y$, an element $x$ of $L$ is extremal
if and only if there is a linear functional $f_x: y\mapsto f_x(y)$
on $L$ such that, for all $y\in L$, we have
\begin{equation}\label{exte}
[x,[x,y]]= f_x(y) x.
\end{equation}
Note that $f_x(y)=0$ if $x$ and $y$ commute. 
We also write $f(x,y)$ instead of $f_x(y)$. 

\begin{Lm}\label{3.1}
If $x, y\in \Ec$ then $f_x(y)=f_y(x)$.
\end{Lm}

\begin{proof}
By \eqref{exte}, $[y,[x,[x,y]]] = f_x(y) [y,x]$, but by Jacobi,
 antisymmetry, and \eqref{exte},
\[[y,[x,[x,y]]] = -[x,[[x,y],y]] -[[x,y],[y,x]]=-[x,[y,[y,x]]]=-f_y(x)[x,y].\]  
This shows that the lemma holds if $x$ and $y$ do not commute.
Otherwise, $f_x(y) = 0 = f_y(x)$.
\end{proof}

\bigskip
We give some general identities and 
properties.  At least the first two rewriting rules
go back to Premet, see \cite{Che}.

\begin{Lm}\label{ide}\label{3.2}
Let $x \in \Ec$. Then, for all $y, z \in L$,
\begin{eqnarray}
\label{ide1}2[[x,y],[x,z]] &=&  f_x([y,z])x + f_x(z)[x,y] - f_x(y)[x,z],\\
\label{ide2}2[x,[y,[x,z]]]&=& f_x([y,z])x - f_x(z)[x,y] - f_x(y)[x,z].
\end{eqnarray}
\end{Lm}

\begin{proof}
Let $y, z \in L$. Then 
\[[[x,y],[x,z]] =
 -[[y,[x,z]],x]  -[[[x,z],x],y] = [x,[y,[x,z]]] + f_x(z) [x,y].\]
Similarly, interchanging $y$ and $z$,
\[[[x,y],[x,z]] =  - [x,[z,[x,y]]] - f_x(y)[x,z].\]
Furthermore, 
\[[x,[z,[x,y]]] = -[x,[x,[y,z]]] - [x,[y,[z,x]]] =  -f_x([y,z]) x + [x,[y,[x,z]]].\]
Hence 
\[ [[x,y],[x,z]] = [x,[y,[x,z]]] + f_x(z)[x,y],\]
\[ [[x,y],[x,z]] = f_x([y,z])x - [x,[y,[x,z]]] - f_x(y) [x,z]. \]

{From} this the two equations follow. 
\end{proof}

\begin{Lm}\label{rel}\label{3.3}
Let $x,y\in \Ec$ and $z \in L$. Then

\begin{equation}\label{rela}
2[[x,y],[x,[y,z]]] = f_y(z)f_x(y)x + f_x([y,z])[x,y] - f_x(y)[x,[y,z]],
\end{equation}

\begin{multline}\label{relb}
2[[x,y],[[x,y],z]] = (f_x([y,z])-f_y([x,z]))[x,y]\\
+f_x(y)\left(f_x(z)y+f_y(z)x-[y,[x,z]]-[x,[y,z]]\right).
\end{multline}
\end{Lm}

\begin{proof} 
Expression \eqref{rela} is the first identity of Lemma \ref{3.2}
with $[y,z]$ replacing $z$.

Identity \eqref{relb} follows from two applications of \eqref{rela} and Jacobi:

\begin{eqnarray*}
2[[x,y],[[x,y],z]]&=&2[[y,x],[y,[x,z]]]+2[[x,y],[x,[y,z]]]\\
&=&f_x(z)f_x(y)y + f_y([x,z])[y,x]-f_x(y)[y,[x,z]]\\
&&\quad +f_y(z)f_x(y)x + f_x([y,z])[x,y]-f_x(y)[x,[y,z]]\\
&=& (f_x([y,z])-f_y([x,z]))[x,y]\\
&&\quad +f_x(y)\left( f_x(z)y+f_y(z)x-[y,[x,z]]-[x,[y,z]]\right).
\end{eqnarray*}
\end{proof}

\bigskip
The lemmas above also make clear why we are assuming that the
characteristic of $k$ is not $2$.

\begin{Cor}\label{3.4}
If $x,y\in \Ec$ with $f_x(y) = 0$ but $[x,y] \neq 0$, then $[x,y] \in \Ec$,
with $f_{[x,y]}(z) = \frac{1}{2}(f_x([y,z])-f_y([x,z]))$ for $z\in L$.
\end{Cor}

\begin{proof} 
Immediate from Identity \eqref{relb} of Lemma \ref{rel}. 
\end{proof}

\bigskip
Let $x\in \Ec$.
Because $\ad_x^3 =0$, the exponential of the derivation $s\ad_x$
for $s\in k$
is given by
\[\exp(x,s) := 1+ s \ad_x+ \frac{s^2}{2}\ad_x^2.\]
This is an automorphism of $L$ in view of Identity \eqref{ide1} of Lemma \ref{3.2}.
Denote by $G$ (or $G(L)$ if necessary) the group of automorphisms 
of $L$ generated by all $\exp(x,s)$ for $x\in \Ec$ and $s\in k$. 

We use these automorphisms to show that if a Lie algebra is generated 
(as a Lie algebra) by extremal elements, it is also linearly 
spanned by extremal elements.  

\begin{Lm}\label{linspan}\label{3.5}
If $L$ is generated as a Lie algebra by extremal elements, then it is linearly
spanned by the set $\Ec$ of all extremal elements.   
\end{Lm}

\begin{proof}
We use induction on the length of $z\in L$ as a bracketing of elements
from $\Ec$ (that is, obtained by forming successive brackets of
elements from $\Ec$).  It readily reduces to the case where this
length is 2, the case of length 1 being $z\in \Ec$.  Suppose $z=[x,y]$
with $x,y \in \Ec$.  Now let $m = \exp(x,1)y = y
+[x,y]+\frac{1}{2}f_x(y)x$.  As $\exp(x,1)$ is an automorphism, $m$ is
extremal and so $[x,y]$ is in the linear span of the extremal elements
$x,m,y$.
\end{proof}

\bigskip
We show how the map $f$ defined in the beginning of this section
gives rise to a symmetric bilinear 
associative form.

\begin{Thm} \label{3.6}
Suppose that $L$ is generated by $\Ec$.  There is a unique
bilinear symmetric form $f: L\times L\to k$ such that, for each $x\in
\Ec$, the linear form $f_x$ coincides with $y\mapsto f(x,y)$. This
form is associative, in the sense that $f(x,[y,z]) = f([x,y],z)$ for
all $x,y,z\in L$.
\end{Thm}

\begin{proof}
Take $ \{ u_i \}_i$ to be a basis of $L$ consisting of elements from $\Ec$.
Its existence is guaranteed by Lemma~\ref{3.5}.
Notice  that if $x$ is extremal, then  
$[\lambda x,[\lambda x,y]]=\lambda f(x,y)\lambda x$ and 
so $\lambda x$ is also extremal with 
$f(\lambda x,y)=\lambda f(x,y)$.  
Now suppose $x=\sum_i\lambda_i u_i$.  Define $f_x$ by 
$\sum_i \lambda_if_{u_i}.$  Notice from above this is 
$\sum_i f_{\lambda_iu_i}$.  Now suppose $\sum_i u_i=\sum_i v_i$ 
are two ways of writing an
element $x$ of $L$ as a sum with $v_1,\ldots,v_m$, 
extremal elements.  By Lemma \ref{3.1}, for $z\in \Ec$, we have 
\[ f_x(z)=\sum_i f_{u_i} (z)= \sum_i f_z (u_i)=f_z(\sum_iu_i)=f_z(\sum_iv_i)=
\sum_i f_z (v_i)=\sum_i f_{v_i} (z).\]
Since $\Ec$ spans $L$, we conclude
$\sum_i f_{u_i}=\sum_i f_{v_i}$. Thus, $f_x$ is a well
defined linear functional.  It readily follows that $f(x,y)=
f_x(y)$ defines a bilinear form.  It is symmetric by Lemmas \ref{3.1}
and \ref{3.5}. 

It remains to establish that $f$ is associative.
Take $x,y,z\in \Ec$.
Interchanging $x$ and $y$ in \eqref{relb} of Lemma \ref{3.3}, we obtain
\[2[[x,y],[y,[x,z]]] = - f(x,z)f(x,y)y + f(y,[x,z])[x,y] + f(x,y)[y,[x,z]].\]
On the other hand, by Jacobi and Lemma \ref{ide},
\begin{eqnarray*}
2[[x,y],[y,[x,z]]] &=& -2 [y,[[x,z],[x,y]]]-2[[x,z],[[x,y],y]]\\
&=& \left[ y, f(x,[y,z])x+f(x,z)[x,y]-f(x,y)[x,z] \right]\\
& & -2f(x,y)[[x,z],y]\\
&=& -f(x,[y,z])[x,y] -f(x,z)f(x,y)y +f(x,y)[y,[x,z]].
\end{eqnarray*}
Suppose now $[x,y]\neq 0$.
Comparing the coefficients of $[x,y]$ in the two expressions for
$2[[x,y],[y,[x,z]]]$, we find $-f(x,[y,z])=f(y,[x,z])$.
But then, by symmetry of $f$, we have
$f(x,[z,y])=-f(x,[y,z])=f(y,[x,z])=f([x,z],y)$, proving 
$f(x,[z,y])=f([x,z],y)$ whenever $[x,y] \neq 0$.
Similarly, $f(x,[y,z]) = f([x,y],z)$, whenever $[x,z] \neq 0$ and
$f(y,[x,z]) = f([y,x],z)$, whenever $[y,z] \neq 0$.

Suppose that $[x,y] = 0$.
We have to show $f(x, [z,y]) = f([x,z],y)$.
If $[x,z] \neq 0$ and $[y,z] \neq 0$, then $f(x,[z,y]) = 0 = f([x,z],y)$
by the previous equations. Hence (up to symmetry) we are left with the case
$[x,y] = 0$, $[x,z] = 0$. Then $f([x,z],y) = 0$, so we must show that
$f(x,[z,y]) = 0$. By Jacobi, $[x,[y,z]] = 0$. Applying $\ad_x$ yields
$f(x,[y,z]) = 0$ as needed.

We conclude that $f(x,[z,y]) = f([x,z],y)$ for all $x,y,z\in \Ec$. Since
$\Ec$ linearly spans $L$, this proves that $f$ is associative.
\end{proof}

\begin{Cor}\label{zero}\label{3.8}
If $x,y\in \Ec$ with $f_x=0$ and $[x,y]\ne0$, then $[x,y] \in \Ec$ with $f_{[x,y]}=0$.
\end {Cor}

\begin{proof}
By Corollary \ref{3.4}, we have $[x,y]\in \Ec$.
Since $f$ is associative by Theorem \ref{3.6},
we find $f_{[x,y]}(z) =
f([x,y], z) = f(x, [y,z]) = 0$ for $z\in L$.
\end{proof}

\section{Some examples}
\label{reisla}


We first discuss Lie algebras generated by two extremal 
elements. 

\begin{Lm}\label{2gen}
Let $L$ be generated by two extremal elements
$x,y\in \Ec$.  Then one of the following three assertions holds.

\begin{itemize}
\item[(i)] $L=kx+ky$ is Abelian, $\Ec = L\setminus\{0\}$, and $f = 0$.
     
\item[(ii)] $L=kx+ky+kz$ with $z=[x,y]\ne0$, its center $Z(L)$
and commutator subalgebra $[L,L]$ coincide with $kz$ 
and $L/Z(L)$ is Abelian. Thus, $L$ is a Heisenberg algebra.
Moreover $\Ec = L\setminus\{0\}$ and $f=0$. 

\item[(iii)] $L\cong \got \sli_2$ and $\Ec$ consists of all nilpotent elements of
$L$.
Thus, $$\Ec \cup \{0\} = kx \cup ky \cup \bigcup_{\delta\in k \setminus\{0\}}
k\left( \delta x + \delta^{-1} \lambda y + [x,y] \right),$$
where $\lambda=\frac{1}{2}f(x,y)\ne0$. 
\end{itemize}
\end{Lm}

\begin{proof}
Clearly, $x$, $y$ and $[x,y]$ linearly span $L$. The rest is straightforward.
\end{proof}

\begin{Example}  
To see that there are real differences in the characteristic 2 case,
let us take $k=\ZZ/2\ZZ$ and look at the case where $L$ is generated by 
two \spe{} elements $x,y$ with $[x,[x,y]]=x$. Then $L=kx+k[x,y]+ky$ is not
isomorphic to $\sli_2$, as it is simple whereas $\sli_2$ is not.
\end{Example}

\bigskip
Next we deal with finite dimensional Lie algebras of Chevalley type.
A long root element of such a Lie algebra is an element of the form
$x_\alpha$ with $\alpha$ a long root, or the image  of such an element
under a Lie algebra automorphism.

\begin{Prop}\label{longchev}
Let $L$ be a Lie algebra of Chevalley type over $k$.  Then
$\Ec$ contains the long root elements of $L$.  In particular, $L$ is
generated by $\Ec$.  Moreover, if $k$ has characteristic distinct from
2 and 3, every element of $\Ec$ is a long root element.
\end{Prop}

\begin{proof}
Consider a Chevalley basis $B$ of $L$ with respect to a given Cartan subalgebra
$H$ of $L$, write $\Phi$ for the corresponding root system,
and $x_\beta$ for the element of $B$ associated with
a given root $\beta\in\Phi$.
Suppose that $\alpha$ is a long root in $\Phi$.
Then, for $h\in H$, we have
$[x_\alpha,[x_\alpha,h]] = \alpha(h) [x_\alpha,x_\alpha] = 0$.
Moreover, for $\beta\in \Phi$,
the sum $\beta + 2 \alpha$ is not a root in $\Phi$ unless
$\beta = -\alpha$, so
$[x_\alpha,[x_\alpha,x_\beta]] = 0$ unless $\beta = -\alpha$, in which
case the result is a multiple of $x_\alpha$.

\bigskip
To see that $L$ is generated by the long root elements, observe that
$L$ is generated by the root elements of the basis $B$ and (by
analysis of the Lie rank 2 case) that any short root element can be
written as a sum of three long root elements, as is clear from the
following computations in the Lie algebras of type $B_2$ and $G_2$,
respectively.

\begin{eqnarray*}
{\rm for}\ B_2:&\quad&
\underbrace{\exp(x_{\eps_1},1)x_{-(\eps_1-\eps_2)}}_{\rm long} =
\underbrace{x_{-(\eps_1-\eps_2)}}_{\rm long} - \underbrace{x_{\eps_2}}_{\rm
short} + \underbrace{x_{\eps_1+\eps_2}}_{\rm long}\\ {\rm for}\
G_2:&\quad&\underbrace{\exp(x_\alpha,1)x_\beta}_{\rm long}
+\underbrace{\exp(x_\alpha,-1)x_\beta}_{\rm long} = 2\underbrace{x_\beta}_{\rm
long} - 2\underbrace{x_{2\alpha+\beta}}_{\rm short}  
\end{eqnarray*}

As for the last assertion, let $e\in \Ec$.  Under the given
characteristic restriction for $k$, the proof of Lemma 2.1  of
\cite{Ben} and
the paragraph following Step 1 in the proof of Theorem 3.2 of
\cite{Ben} show that there exists a Cartan subalgebra of $L$ with
eigenspace $ke$.  Thus $e$ is a root element.  If the diagram of $L$
is not simply laced, short root elements of $L$ are easily seen not to
be extremal. Hence $e$ is a long root element.
\end{proof}

\bigskip
We expect that, in characteristic 3,
each extremal element is also a long root element.
Such an assertion might follow from a classification of nilpotent elements
in the Lie algebras of Chevalley type, 
but we have not been able to verify this completely.

\section{$L$ is finite dimensional in the finitely generated case}
\label{fingen}
\label{finsec}

We use the work of Zel'manov and Kostrikin \cite {ZK} to show
the following main result.

\begin{Thm}\label{main}
If $L$ is generated as a Lie algebra by a finite number of extremal elements,
then  $L$ is finite dimensional. 
\end{Thm}

\begin{proof} 
Suppose that $L$ is generated by extremal elements $x_1, \dots, x_r$. Denote
by $f$ the associated bilinear form. 
We first consider the case in which $f$ is identically  zero.

\begin{Lm} \label{finite1}\label{4.2}
Suppose that $L$ is generated by elements  $x_1, \dots, x_r$ with $(\ad_{x_i})^2 = 0$ 
for $i=1, \dots, r$.
Then $L$ is finite dimensional and is nilpotent.  Furthermore, there is a
universal such Lie algebra, $\Lc_r$.
\end{Lm}   

\begin{proof} 
Observe that each $x_i$ is in $\Ec$.
By Corollary \ref{zero} $(\ad_w)^2 = 0$ for any bracketing $w$ in
$x_1, \dots, x_r$. (We consider also $x_i$ as a bracketing.)  Since
$\ch k \neq 2$, this implies that $\ad_w \ad_u \ad_w = 0$ for any two
of these bracketings, see \eqref{ide2} of Lemma \ref{3.2}.  In the language of
Zel'manov and Kostrikin \cite{ZK}, each $\ad_w$ is a crust of a thin
sandwich, and by \cite[Theorem 2]{ZK} the associative algebra
generated by $\ad_{x_1}, \dots, \ad_{x_r}$ is nilpotent.  Hence all
products of some fixed length, say $l$, of the $\ad_{x_i}$ are $0$.
This yields that $L$ is a nilpotent Lie algebra with bracketings of
length $l+1$ being $0$. In particular, $L$ is finite dimensional.

To see that there is a universal such Lie algebra, let $F$ be the free 
Lie algebra generated by $r$ elements $f_1,f_2,\dots,f_r$.  Now
factor out by the ideal $J$ generated by all $[f_i,[f_i,u]]$ where
$u$ is a bracketing in the $f_i$.  Now $F/J$ is a universal Lie
algebra spanned by $r$ extremal elements with $f$ identically $0$.  
This completes the proof of the lemma.  
\end{proof}

\bigskip
We now show how this implies Theorem \ref{main}.

\begin{Lm}\label{4.3}
Suppose that $L$ is generated by $r$ extremal elements 
where the values of $f$ need not all be $0$.  Then the dimension of $L$ is at
most $\dim \Lc_r$.  In particular $L$ is finite dimensional.
\end{Lm}

\begin{proof}  
For a nonzero bracketing in $F$, its length is understood to be the
total number of the $f_i$ which appear, counting multiplicities (e.g.,
the length of $[[f_1, f_2],[f_1,f_3]]$ is 4).  Choose a basis $\{w_1 +
J, \dots, w_t + J\}$ of $F/J$ and let $U$ be the linear span of $w_1,
\dots, w_t$. Then $F = U + J$.  Consider the surjective homomorphism
$\wt{\phantom{a}}: F \to L$ determined by $f_i \mapsto x_i$. We claim
that $\{\wt{w_1}, \dots, \wt{w_t}\}$ linearly spans $L$, i.e., $\wt{U}
= L$.

Any element of $J$ is a linear combination of terms
\begin{eqnarray}\label{ffu}
[u_s, [ \dots, [u_1, [f_i, [f_i, u]]] \dots ]]
\end{eqnarray}
with $u, u_1, \dots, u_s$ bracketings in the $f_i$.
Suppose that $\wt{U}$ is properly contained in $L$. Then there exists
a nonzero bracketing $ g \in F$ of minimal length with $\wt{g} \not \in \wt{U}$.
We may write $g$ as a linear combination of the $w_i$ and terms $T_i$ 
as in \eqref{ffu} with all $w_i, T_i$ of the same length as $g$. (Note that in the
free Lie algebra $F$, the linear span of the bracketings of a fixed length $\ell$
trivially intersects the linear span of all bracketings of lengths distinct from $\ell$.)
Among these terms $T_i$ there is a term $T =
[u_s, [ \dots, [u_1, [f_i, [f_i, u]]] \dots ]]$ as in \eqref{ffu}
with $\wt{T} \not\in \wt{U}$.
Set $T_0 := [u_s, [ \dots, [u_1, f_i] \dots]]$. Then 
$f(x_i, \wt{u})\wt{T_0} = \wt{T} \not \in \wt{U}$. 
Hence $f(x_i, \wt{u}) \neq 0$ and $\wt{T_0} \not\in \wt{U}$ with $T_0$
a bracketing of 
shorter length than $g$, a contradiction with the choice of $g$.
\end{proof}

\bigskip
Combining these two lemmas proves Theorem \ref{main}.
\end{proof}

\bigskip

\begin{Remark}\label{4.4}
We have worked out the Lie algebra for up to five 
generators and summarize the dimensions below.  Let 
$\Lc_r$ be the universal Lie algebra generated by $r$ extremal 
generators subject to $f$ being identically $0$.  Also,
let $\Rc_r$ be the free associative algebra over
$k$ generated by $r$ elements 
$y_1,y_2,\dots ,y_r$ for which $y_i^2=0$ and $y_iwy_i=0$ 
where $w$ is any bracketing in the $y_i$s. Then $\Lc_r/Z(\Lc_r)$,
the image of $\Lc_r$ under $\ad$, is a quotient of
$\Rc_r$ in view of \eqref{ide2} of Lemma \ref{3.2}. Our computations 
have found the values of the corresponding dimensions in the 
table below.  
\end{Remark}

{\small \begin{center}
\bigskip
\begin{tabular}{|c|c|c|}
\hline
Number of generators& $\dim  \Lc_r$ &$\dim \Rc_r$ \\ 
\hline\hline
1 &1 &2\\ \hline
2 &3 & 5\\ \hline
3 &8 &19\\ \hline
4 &28 &193\\ \hline
5 &537 & ?\\ \hline
\end{tabular}
\end{center}}

\bigskip

By way of example, consider the case $r=2$. Put $A_1=\ad_{x_1}$ and
$A_2=\ad_{x_2}$.  Then a basis for $\Lc_2$ is $x_1,x_2,[x_1,x_2]$ of
size 3.  A basis for $\Rc_2$ is $I$, $A_1$, $A_2$, $A_1 A_2$, $A_2A_1$,
of size 5, as can be seen immediately. This explains the corresponding entries
for $r=2$ in the above table.

We will explain in later sections how some of the other entries have
been determined.  For the remainder we refer to the algorithmic
methods described in \cite{Reut} and \cite{LR} using Lyndon words, and
to the computational algebra packages GAP \cite{gap} (including the FPSLA program by
Gerdt \& Kornyak), and {\sc LiE}
\cite{Lie}.

\bigskip
Let 
$\Rad(f)$ denote the radical of $f$, that is, the set of all $y \in L$ for which
$f(y,z)=0$ for all $z\in L$. 
\bigskip

\begin{Cor}
If $L$ is simple, then $f$ is non-degenerate and there is no $x \in L$, 
$x \neq 0$, with $\ad_x^2 = 0$.
(This means that $L$ is non-degenerate in the sense of \cite{Ben}.)
\end{Cor}

\begin{proof}
We know that
$f$ is non-trivial by Lemma~\ref{4.2}. 
(Otherwise $f$ is identically $0$ and $L$ is nilpotent.)
Thus, $\Rad(f)$ is a proper ideal of $L$
and so, as $L$ is simple, $\Rad(f) = 0$. This means that $f$ is non-degenerate.

Suppose $x \in L$ with $\ad_x^2 = 0$. Then
$x$ is extremal and $[x,[x,y]] = 0$ for all $y \in L$.
Hence $x \in \Rad(f)$. This means $x = 0$. 
\end{proof}

\begin{Remark}
By the above corollary and the proof of Theorem 3.2 in \cite{Ben}, for
$p>5$, the only modular finite dimensional simple Lie algebras over an
algebraically closed field of characteristic $p$ and generated by
extremal elements, are quotients of those of Chevalley type.
\end{Remark}          

\section{The 3 generator case}
\label{3gensec}

Suppose that $L$ is generated by three extremal elements $x$, $y$, $z$. 
We wish to determine the various
possibilities for $L$.
The identities of Lemmas \ref{ide} and \ref{rel},
together with the identity
\begin{multline*}
2[[x,[y,z]],[y,[x,z]]] = -\frac{1}{2} \left
( f_y(z)f_x([y,z])x+f_x([y,z])f_x(z)y+f_x([y,z])f_x(y)z \right)\\ 
- f_y(z)f_x(z)[x,y]  + f_y(z)f_x(y)[x,z] - f_x(z)f_x(y)[y,z],
\end{multline*}
show that
\[x,y,z,[x,y],[x,z],[y,z],[x,[y,z]],[y,[x,z]]\]
linearly span $L$. In particular, $L$ is at most 8-dimensional.
Hence, $\Lc_3$ 
is also of dimension at most $8$.
It is readily checked though,
that the above 8 bracketings provide 
a basis in the free case, so $\dim \Lc_3 = 8$.

\begin{Example}
\label{sl3ex}
The algebra $\sli_3$ can be generated by $3$ elements.  
It can be realized as follows:
$$x = \left( \begin{array}{ccc} 0&1&0\\ 0&0&0\\ 0&0&0 \end{array} \right), 
y = \left( \begin{array}{ccc} 0&0&0\cr 1&0&0\cr 0&0&0 \end{array} \right),
z = \left( \begin{array}{ccc} 1&1&1\cr 1&1&1\cr -2&-2&-2 \end{array} \right).$$
In this example, we have 
$$f(x,y)=f(x,z)=f(y,z) =-2\quad\hbox{and}\quad f(x,[y,z])=0.$$
\end{Example}

\bigskip
The actions of $\ad_ x$ and $\ad_ y$ on the linear generators of $L$ can
be fully described by means of the identities in terms of the four
parameters $ f(x,y)$, $ f(x,z)$, $ f(y,z)$, and $ f(x,[y,z])$.

\bigskip
We can describe the four parameters pictorially by drawing a triangle
with vertices $x$, $y$, $z$, and labeling the edge $\{x,y\}$ by the
`edge' parameter $f(x,y)$ and so on, and putting the `central'
parameter $f(x,[y,z])$ in the middle, with an indication of
orientation (note that $f(x,[y,z])=f([x,y],z)=f(z,[x,y])$, so the value
is invariant under cyclic permutations of the nodes).

\bigskip
We shall reduce $f(x,[y,z])$ to zero by transforming the generators using
elementary transformations. 
To begin reduction, consider the triple $x,y,\exp(x,s)z$, where $s\in
k$.  It has parameters:
\begin{eqnarray*}
f(x,y)&=&f(x,y) \\ 
f(x,\exp(x,s)z)&=&f(x,z)\\
f(y,\exp(x,s)z)& =& f(y,z)-sf(x,[y,z])+\frac{1}{2}s^2f(x,y)f(x,z)\\
f(x,[y,\exp(x,s)z]) &=& f(x,[y,z])-sf(x,z)f(x,y).
\end{eqnarray*}
Clearly, this triple again consists of 
\spe{} elements, and generates the same algebra as $x$, $y$ and $z$.
If at least two of the three edges have nonzero labels
(e.g., $f(x,z)$ and $f(x,y)$), then we can transform the central parameter
to $0$ (by taking $s = f(x,[y,z])/(f(x,z)f(x,y))$).

On the other hand, if at most one edge is nonzero, say $f(x,y)$,
and the central parameter $f(x,[y,z])$ is also nonzero,
then the above transformation shows that we can move
to three \spe{} generators $x$, $y$, $\exp(x,s)z$ with
one more edge (namely $f(y,\exp(x,s)z)$ and, if applicable, $f(x,y)$) nonzero. Hence,
we can reduce to the previous case (if necessary in two steps),
and so we may assume that the central parameter is zero: $f(x,[y,z])=0$. 

Next, we scale $x$, $y$, $z$ to $\alp x$, $\beta y$, $\gam z$
for nonzero $\alp,\beta,\gam\in k$.
This leaves $f(\alp x, \beta  \gam[y, z])=0$ and
changes the edge labels to
$$\alp\beta f(x,y),\ \alp\gam f(x,z),\ \beta\gam f(y,z).$$
We claim that, at the cost of a field extension of $k$,
all nonzero edge labels may be transformed into $-2$.
If at least one of them is zero, this is obvious.
Otherwise, take 
$$\alp = \sqrt{ \frac{-2f(y,z)}{f(x,y)f(x,z)} }, \enspace
 \beta = \sqrt{ \frac{-2f(x,z)}{f(x,y)f(y,z)} }, \enspace
  \gam = \sqrt{ \frac{-2f(x,y)}{f(x,z)f(y,z)} }.$$

Thus we are left with four essentially different cases, distinguished
by the number of nonzero labeled edges in the triangle.
Straightforward computation using GAP \cite{gap} leads to the following
descriptions of the resulting four  Lie algebras.

\begin{Thm}\label{3gen}
Suppose that $L$ is generated by three \spe{} elements. Then after
extending the field if necessary, $L$ is generated by three \spe{}
elements whose central parameter is zero, and whose nonzero edge
parameters are $-2$. In particular, $L$ is a quotient of a Lie algebra
$M$ generated by extremal elements $x$, $y$, $z$ with $f(x,[y,z]) = 0$
and $\dim M = 8$. Moreover, according as the number of nonzero edge
parameters is $0$, $1$, $2$, $3$, the Lie algebra $M$ has the
following form.

\begin{itemize}
\item[(0)] $f=0$ and $M \cong \Lc_3$. Thus,
$M$ is nilpotent, with $[M,M]=k[x,y]+k[x,z]+k[y,z]+Z$ where
  $Z=k[x,[y,z]]+k[y,[x,z]]=[[M,M],M] $ is the center of $M$. 

\item[(1)] $f(x,y)=-2$, $f(x,z)=f(y,z)=0$,
and $M = Z \oplus [M,M]$ where  $Z=k(z-[x,[y,z]]-[y,[x,z]])$
is the center of $M$. The solvable radical of
  $M$ is $R =kz+k[x,z]+k[y,z]+k[x,[y,z]]+k[y,[x,z]]$. The subalgebra
  $S=kx+k[x,y]+ky$ is isomorphic to $\sli_2$ and $M$ is the
  semi-direct product of $S$ and $R$. The $S$-modules
  $k[x,z]+k[y,[x,z]]$ and $k[y,z]+k[x,[y,z]]$ are irreducible.

\item[(2)] $f(x,y)=f(x,z)=-2$, $f(y,z)=0$,
and  $M$ is the semi-direct product of
  $S=kx+k[x,y]+ky$, which is isomorphic to $\sli_2$,
and the solvable radical $R$ of
  $M$. Moreover, $R =k(y-\frac{1}{2}[y,[x,z]]) + k(z-\frac{1}{2}[y,[x,z]] ) +
  k([x,y]-[x,z]) + k[y,z] + k[x,[y,z]]$, 
$[R,R]=k(y+z-[y,[x,z]])+k[y,z]+k[x,[y,z]]$ and $[R,[R,R]] =
  k[y,z]+k[x,[y,z]]$. 
The center of $M$ is trivial and $[M,M] =
  M$. The subspace $k(y+z-[x,[y,z]]-[y,[x,z]])$ of $ [R,R]$ is
 centralized by $S$. 

\item[(3)]  $f(x,y)=f(x,z)=f(y,z)=-2$
and $M \cong \sli_3 $ as described in Example \ref{sl3ex}.
\end{itemize}
\end{Thm}

\bigskip
The algebra $\Rc_3$ can be determined easily.  Again let 
$A_i=\ad_{x_i}$.  
\bigskip

{\small
\begin{center}
\begin{tabular}{|c|c|c|} \hline 
Words & Conditions &Number \\ 
\hline \hline
$I$ &identity &1\\ \hline
$A_i$ & & 3\\ \hline
$A_i A_j$ &$i$, $j$ distinct  &6\\ \hline
$A_i A_j A_k$ &$i$, $j$, $k$ distinct &6\\ \hline
$A_i A_j A_k A_i$ &$i$, $j$, $k$ distinct & 3\\ \hline
\end{tabular}
\end{center}}

\bigskip
Notice that in the last line as $A_i[A_j,A_k]A_i=0$, we get 
$A_iA_jA_kA_i=A_iA_kA_jA_i$ and so there are only $3$ of these.  Notice 
these are the only possibilities as multiplying one of the words of 
length $4$ by any $A_l$ obviously gives $0$ on one side and gives 
$0$ on the other side after using the commutation rule just given.  
There is such an algebra, as we found using the $4$ generator 
Lie algebra $\Lc_4$ to be described in the next section.   

\section{The 4 generator case}
\label{4gensec}

\begin{defn}
A {\em monomial of length $s$} is a bracketing of the form
\[ [x_1,[x_2, \dots [x_{s-1},x_s] \dots ]].\]
A monomial is {\em reducible} if it is a linear combination of monomials of
strictly smaller length.
\end{defn}

\begin{Lm}
Let $L$ be generated by a subset $D$.
Then $L$ is the linear span of all monomials in elements of $D$.
\end{Lm}

\begin{proof}
By induction on the length of a bracketing (with respect to $D$) and Jacobi.
\end{proof}

\begin{Prop}\label{6.3}
Let $L$ be generated by the extremal elements $x, y, z, u$.
Then $L$ is linearly spanned by the following $28$ monomials of length
$\leq 5$:
\[x, y, z, u,\]
\[[x,y], [x,z], [x,u], [y,z], [y,u], [z,u], \]
\[[x,[y,z]], [x,[y,u]], [x,[z,u]], [y,[x,z]], [y,[x,u]], [y,[z,u]], [z,[x,u]],
[z,[y,u]], \]   
\[[x,[y,[z,u]]], [x,[z,[y,u]]], [y,[x,[z,u]]], [y,[z,[x,u]]], [z,[x,[y,u]]],
[z,[y,[x,u]]], \]  
\[[x,[y,[z,[x,u]]]], [y,[x,[z,[y,u]]]], [z,[x,[y,[z,u]]]],
[u,[x,[y,[z,u]]]]. \] 
\end{Prop}

\begin{proof}
Since $L$ is linearly spanned by the monomials in $x, y, z, u$, we have to
show that each monomial may be written as a linear combination of the
given 28 elements.

All monomials of length 1 are on the list. There are $4\cdot 3$ monomials
of length 2 with different factors. 
Since $[a,b] = -[b,a]$ for all $a, b \in L$, all of them may be
expressed by the 6 monomials of length 2 on the list.

All monomials of length 3 which involve only 2 letters are reducible, since
$x, y, z, u$ are extremal.
There are $4\cdot 3 \cdot 2$ monomials of length 3 with three different
letters. With antisymmetry, all $[x,[a,b]]$ and all $[y,[a,b]]$ may be
expressed. By  Jacobi this holds also for the remaining ones.

All monomials of length 4 which involve only three letters are
reducible, see the identity for $2[x,[y,[x,z]]]$ in Lemma~\ref{ide}.
There are $4 \cdot 3 \cdot 2$ monomials of length 4 with different factors.

With Jacobi we have the following equations
\begin{equation}\label{(1)}
[x,[y,[z,u]]] - [x,[z,[y,u]]] + [x,[u,[y,z]]] = 0, 
\end{equation}
\begin{equation}\label{(2)}
[y,[x,[z,u]]] - [y,[z,[x,u]]] + [y,[u,[x,z]]] = 0,
\end{equation}
\begin{equation}\label{(3)}
[z,[x,[y,u]]] - [z,[y,[x,u]]] + [z,[u,[x,y]]] = 0.
\end{equation}

The elements in the first two columns are on the list.
Hence all monomials beginning with $x$, $y$, or $z$ may be expressed
as a linear combination of the given 28 basis elements.
For the monomials beginning with $u$, we calculate:
\begin{equation}\label{(4)}
[u,[x,[y,z]]] - [x,[u,[y,z]]] + [[y,z],[u,x]] = 0,
\end{equation}
\begin{equation}\label{(5)}
[u,[y,[x,z]]] - [y,[u,[x,z]]] + [[x,z],[u,y]] = 0,
\end{equation}
\begin{equation}\label{(6)}
[u,[z,[x,y]]] - [z,[u,[x,y]]] + [[x,y],[u,z]] = 0.
\end{equation}

The products of the form $[[a,b],[c,d]]$ may be expressed as follows.
\begin{equation}\label{(7)}
[x,[y,[z,u]]] - [y,[x,[z,u]]] + [[z,u],[x,y]] = 0,
\end{equation}
\begin{equation}\label{(8)}
[x,[z,[y,u]]] - [z,[x,[y,u]]] + [[y,u],[x,z]] = 0,
\end{equation}
\begin{equation}\label{(9)}
[y,[z,[x,u]]] - [z,[y,[x,u]]] + [[x,u],[y,z]]= 0.
\end{equation}

Hence also the monomials of length $4$ beginning with $u$ may be expressed
as a linear combination of the monomials on the list.

We obtain monomials of length 5 by multiplying a letter from the left to a
monomial $[a,[b,[c,d]]]$ of length 4. This yields four possibilities: First
$[a,[a,[b,[c,d]]]]$, which obviously is reducible. Second $[b,[a,[b,[c,d]]]]$,
which is of the form \linebreak $[b,[a,[b,e]]]$ and hence it is
reducible. Third $[c,[a,[b,[c,d]]]]$. There is no obvious way to rewrite
this. Fourth $[d,[a,[b,[c,d]]]]$, which yields the previous case by
interchanging $c$ and $d$. We are left with the monomials of the form
$[c,[a,[b,[c,d]]]]$ or $[d,[a,[b,[c,d]]]]$, where $[a,[b,[c,d]]]$ is one of
the monomials of length 4 on the list. This yields the following 12 monomials
$m_{11}, m_{12}, \dots, m_{34}$ ($m_{ij}$ is in row $i$ and column $j$): 
\begin{multline*}
[z,[x,[y,[z,u]]]], [y,[x,[z,[y,u]]]], [z,[y,[x,[z,u]]]], [x,[y,[z,[x,u]]]],\\
[y,[z,[x,[y,u]]]], [x,[z,[y,[x,u]]]],[u,[x,[y,[z,u]]]], [u,[x,[z,[y,u]]]], \\
[u,[y,[x,[z,u]]]], [u,[y,[z,[x,u]]]], [u,[z,[x,[y,u]]]], [u,[z,[y,[x,u]]]].
\end{multline*}

Our intended basis vectors are $m_{11}, m_{12}, m_{14}, m_{23}$.
We may express the other 8 elements as follows:

\smallskip
\begin{tabular}{ll}
$m_{24}$: & with \eqref{(1)},\\
$m_{13}, m_{31}$: & with \eqref{(7)} and \eqref{(6)},\\
$m_{32}$: & with \eqref{(2)},\\
$m_{21}, m_{33}$: & with \eqref{(8)} and \eqref{(5)},\\
$m_{22}, m_{34}$: & with \eqref{(9)} and \eqref{(4)}. 
\end{tabular}

\medskip
Finally, we show that all monomials of length 6 are reducible.
We multiply the monomials of length 5 on the list with a letter.
Note that all these monomials of length 5 are of the form
$\pm [c,[a,[b,[c,d]]]]$. Multiplication from the left with $c$ or $a$ yields
reducible monomials.  
Next, we deal with $[b,[c,[a,[b,[c,d]]]]]$.
With Jacobi we may pass from $[a,[b,[c,d]]]$ to $[b,[[c,d],a]]$ and $[[c,d],[a,b]]$. {From}
$[[c,d],[a,b]]$, we pass to $[d,[c,[a,b]]]$ and $[c,[[a,b],d]]$. Now we multiply first with
$c$, then with $b$ from the left. This yields products of the form
$[b,[c,[b,[[c,d],a]]]]$, $[b,[c,[d,[c,[a,b]]]]]$ and $[b,[c,[c,[[a,b],d]]]]$, which are all
reducible. (Look for patterns of the type $[u,[v,[u,w]]]$.)
The last monomial we have to reduce is of the form
$[d,[c,[a,[b,[c,d]]]]]$. We pass from $[b,[c,d]]$ by Jacobi to $v_1 = [c,[d,b]]$ and 
$v_2 = [d,[c,b]]$. The products $[d,[c,[a,v_1]]]$ and $[d,[c,[a,v_2]]]$ are reducible.
(In the second one, the two letters $d$ are at distance $3$.)

As a consequence all monomials of length $\geq 6$ are reducible and may be
written as a linear combination of the 28 vectors on the list.
\end{proof}

\begin{Remark}
Removing all elements involving a letter $u$ from the list for the
4-generator case, yields a spanning set of size 8 for the 3-generator
case.  More generally, for $r\in\NN$, $r>1$, the Lie algebra
$\Lc_{r-1}$ is a Lie subalgebra of $\Lc_r$, generated by the first
$r-1$ elements of a set $\{x_1,\ldots, x_r\}$ of extremal generators
for $\Lc_r$.  The module generated by $x_r$ under the ad action of
$\Lc_{r-1}$ gives the algebra $\Rc_{r-1}$ defined in Remark
\ref{4.4}. (This can be proved by observing that the defining
relations for $\Lc_r$ are homogeneous with respect to the multidegree
counting the number of each $x_i$, so that a non-trivial relation
amongst bracketings in $\Rc_{r-1}x_r$ would lead to a non-trivial
homogeneous relation of $x_r$-degree 1. Thus, no relation in
$\Rc_{r-1}x_r$ involves relations coming from $\ad_{x_r}^2 = 0$, and
so any linear relation amongst bracketings in $\Rc_{r-1}x_r$ is a
consequence of the defining relations of $\Rc_{r-1}$.) 

Counting the monomials with a single $u$
in the list of Proposition \ref{6.3}, we find $\dim\Rc_3 = 19$,
distributed according to length as follows, where length of course
stands for the length in $\Rc_3$, which is one less than
the length of these monomials of Proposition \ref{6.3}.

\medskip
{\small
\begin{center}
\begin{tabular}{|c|c|c|c|c|c|} \hline 
\text{length} &0 &1 &2 &3 &4   \\ \hline
\text {number}&1 &3 &6 &6 &3   \\ \hline
\end{tabular}
\end{center}}

\medskip
The algebra $\Rc_4$ has dimension $193$.  This can be seen by 
routine enumeration which we do not include.  However we 
list the numbers of words of each length in the table below.  

{\small
\begin{center}
\begin{tabular}{|c|c|c|c|c|c|c|c|c|c|c|} \hline 
\text{length} &0 &1 &2 &3 &4 &5 & 6& 7& 8& 9  \\ \hline
\text {number} &1 &4 &12 &24 &36 &40 &36 &24 &12 &4  \\ \hline
\end{tabular}
\end{center}}

\medskip
Notice that in all of these algebras $\Rc_2$, $\Rc_3$, $\Rc_4$ the number 
of words with a given length is symmetric about the highest 
number after deleting the one word of length~$0$ which is the
identity.  
We wonder if this behavior persists for $\Rc_m$ with $m\ge 5$.
\end{Remark}

\section{Quadratic modules}
\label{quadsec}

\begin{defn}
Let $L$ be a Lie algebra generated by extremal elements.
We call an $L$-module $U$ quadratic if there is a generating set $D$ of extremal elements of $L$
such that $x \cdot (x \cdot U) = 0$ for all $x \in D$.
\end{defn}

\begin{Remark}
Suppose that $L$ is generated by $n$ extremal
elements. Consider the subalgebra $A$ generated by the set $D$ of the first $n-1$
extremal generators of $L$. Since $A$ is an $A$-submodule of $L$, the
quotient $L/A$ is also an $A$-module.  Clearly, $x \cdot (x\cdot L/A)
= 0 $ for all $x\in D$, and $D\subseteq \Ec(L)\cap A \subseteq \Ec(A)$.
Therefore, $L/A$ is a quadratic $A$-module.
\end{Remark}

\bigskip
We give the irreducible quadratic modules
for Lie algebras of Chevalley type.

\begin{Prop}
Let $L$ be a Lie algebra of Chevalley type and $k$
a field of characteristic distinct from 2 and 3.
The highest weights of its non-trivial quadratic
highest weight modules of finite dimension are given in the table below.
{\small
\begin{center}
\medskip
\begin{tabular}{|c|c|}
\hline
{\rm Type of} $L$ & {\rm Highest weights}\\
\hline\hline
$A_n$ &$\omega_1,\ldots,\omega_n$\\
$B_n$ &$\omega_1$, $\omega_n$\\
$C_n$ &$\omega_1,\ldots,\omega_n$\\
$D_n$ &$\omega_1,\omega_{n-1},\omega_{n}$\\
$E_6$&$\omega_1,\omega_6$\\
$E_7$&$\omega_7$\\
$F_4$ &$\omega_4$\\
$G_2$ &$\omega_1$\\ \hline
\end{tabular}
\end{center}}
\end{Prop}

\begin{proof}{(Sketch;
a strongly related result for groups can be found in \cite{PS}.)}
Let $V$ be a non-trivial quadratic highest weight module for $L$ of
finite dimension.  By Proposition \ref{longchev}, there is a single
$G$-orbit of extremal elements in $L$. Moreover, all these (long root)
elements have the same nilpotency index on $V$ (because the highest
weight representation is equivalent to a composition of itself with
conjugation by an element of $G$).  Consequently, each element $x\in
\Ec(A)$ satisfies $x \cdot (x\cdot V) = 0 $.

In particular, for every long root $\alpha$, the corresponding
Chevalley basis element $x_\alpha$ is nilpotent of index 2. It follows
that if $\mu$ is a weight of $V$, then $\mu-2\alpha$ is not.  Since the set
of weights for $V$ is a convex subset of the coset of the root lattice
containing the highest weight, only small weights can occur. Thus, the
result can be readily established by use of LiE, cf.~\cite{Lie}, and
the following argument for the case of a single root length.
\end{proof}

\begin{Lm}
Suppose that $L$ is a Lie algebra of Chevalley type with only
one root length.  If $U$ is an irreducible quadratic
finite dimensional $L$-module, then $U$ is a minuscule weight
representation.
\end{Lm}

\begin{proof}
Let $\lambda$ be the highest weight of $U$. Suppose $\mu$ is a weight
for $U$. Then there is a path of weights from $\mu$ to $\lambda$, such
that, for each adjacent pair $(\mu_1,\mu_2)$ from the path, the
difference $\mu_2 - \mu_1$ is a positive root. Since there is only one
root length, each of these differences is a root whose root element is
extremal. Now take the fundamental $\sli_2$-triplet containing this
root element. Quadraticity means that the
non-trivial irreducible sub-representations in $U$ of
the corresponding subalgebra isomorphic to $\sli_2$ all have dimension
$2$. But $\mu_1$, $\mu_2$ are in the same representation, so they
belong to a $2$-dimensional module. But then, they are conjugate by an
element of the corresponding subgroup of type $A_1$. This establishes
that $U$ has a basis of weight spaces which are all conjugate, which
can be taken as a definition of minuscule.
\end{proof}

\begin{Remark}
For the Lie algebra of type $E_8$, there are no non-trivial
finite-dimensional quadratic modules. Thus, this Lie algebra can only occur
as a direct component of a bigger Lie
algebra generated by extremal elements.
\end{Remark}

\section{The minimal number of extremal generators for Lie algebras of 
Chevalley type}
\label{mingensec}

Let $\g$ be a Lie algebra of Chevalley type over $k$ (of
characteristic distinct from 2).  As we have seen in Proposition
\ref{longchev}, the long root elements of $\g$ are extremal elements
and $\g$ is generated by these extremal elements.  Write $t(\g)$ for
the minimal number of these extremal generators of $\g$. In this
section, we determine $t(\g)$.  We discuss implications for the group
analog in \S\ref{rootgrpsec}.

\medskip
Fix a Cartan subalgebra $\hk$ of $\g$ and let $\Phi$ be the corresponding root
system. Throughout this section $x_\alpha$ and $h_{\alpha}$ denote,
respectively, the root element and the element in $\hk$ of $\g$ corresponding
to the root $\alpha$ in $\Phi$. Also $\exp(x_\alpha)$ denotes $\exp(x_\alpha,
1)$ in the notation of Section \ref{3}.  
The signs appearing here in the multiplication rules for
the $x_\alpha$ are the ones
arising from the implementations in the relevant packages.  In 
general, the signs depend on the choice of a Chevalley basis.

\begin{Lm}\label{gen}\label{8.1}
Let $\g$ be a Lie algebra of Chevalley type and let $ \alpha_1, \dots,
\alpha_n$ be the simple roots in $\Phi$. Denote by $\beta$ the
root of highest height. Then $\g$ is generated by the root elements
$x_{\alpha_1}, x_{\alpha_2}, \dots, x_{\alpha_n}, x_{-\beta}$.
\end{Lm}

\begin{proof}
The simple
root elements with respect to $\Phi$
generate the subalgebra $\nk = \sum_{\alpha \in
  \Phi^+} \g_{\alpha}$ where $\Phi^+$ is the set of positive roots and
  $\g_{\alpha}$ is a root space. The lemma follows from the fact that the
  $\nk$-submodule of $\g$ generated by the root element corresponding to the
  root of lowest height, coincides with $\g$.  
\end{proof}

\begin{Thm} \label{mingent}
Let $\g$ be a Lie algebra of Chevalley type over the
field $k$ of characteristic distinct from 2. Then the number
$t(\g)$ is as given in the table below.
\end{Thm}

{\small
\begin{center}
\medskip
\begin{tabular}{|c|c|c|}
\hline
Type of $\g$ & $t(\g)$ & condition\\
\hline\hline
$A_n$&$n+1$&$n\ge1$\\
$B_n$&$n+1$&$n\ge3$\\
$C_n$&$2n$&$n\ge2$\\
$D_n$&$n$&$n\ge4$\\
$E_n$&$5$&$n=6,7,8$\\
$F_4$&$5$&\\
$G_2$&$4$&\\
\hline
\end{tabular}
\end{center}}

\bigskip
The proof of the theorem will be given in the rest of this section.

\bigskip

\begin{lemma}   \label{lb}
If $\g$ has an irreducible representation of dimension $N$ and its extremal
elements have rank $m$ in this representation, then the number of extremal
elements generating $\g$  is at least $N/m$.
\end{lemma}

\begin{proof}
If $\g$ is generated by the extremal elements $x_1,\ldots, x_t$, then the
image of $\g$ in $V$, the underlying $N$-dimensional vector space, is
generated by  
\[ \langle \im x_1, \ldots, \im x_t\rangle = \langle \im x_1\rangle + \cdots +
\langle\im x_t\rangle ,\]
of dimension at most $tm$. Thus, irreducibility of the representation implies
$t m \ge N$, whence the lemma.
\end{proof}

\bigskip
Here is an upper bound for groups without multiple bonds.

\begin{lemma}\label{ub}
If the root system of $\g$ has just one root length, then $t(\g)\le n+1$,
where $n$ is the rank of $\g$.
\end{lemma}

\begin{proof} 
This follows immediately from Lemma \ref{gen}.
\end{proof}

\begin{lemma}\label{lowb}
If $\g$ is generated by $t$ extremal elements, then $t(\g) \ge 4$ if $\dim \g \ge
9$, and $t(\g) \ge 5$ if $\dim \g \ge 29$.
\end{lemma}

\begin{proof} 
This is immediate from the dimensions given in Section \ref{3gensec}
and in Proposition \ref{6.3}.
\end{proof}

\bigskip
The previous lemmas suffice for the proof of all exact lower bounds on $t(\g)$.

For $A_n$, consider the natural representation of dimension $N = n+1$.
Since extremal elements have rank 1 in this representation,
Lemma \ref{lb} gives $t(A_n) \ge n+1$. 

\medskip
For $B_n$, Lemma \ref{lb} applied to the natural module for $\g$ (of dimension
$2n+1$, in which the extremal elements have rank 2), we find $t(\g) \ge
(2n+1)/2$, whence $t(B_n) \ge n+1$. Observe that for $n=2$, this bound is not
sharp. Since the Lie algebra of type $B_2$ has dimension 10, at least 4
extremal generators are needed in view of Lemma \ref{lowb}. This lower bound
coincides with the result for the Lie algebra of type $C_2$ (which is
isomorphic to the one of type $B_2$). 

\medskip
For $C_n$, Lemma \ref{lb} applied to the natural module for $\g$ (of dimension
$2n$, in which the extremal elements have rank 1), we find $t(C_n) \ge 2n$.

\medskip
For $D_n$, recall that the Lie algebra has a natural representation of
dimension $2n$, in which extremal elements have rank 2; hence, by Lemma
\ref{lb}, $t(D_n) \ge n$. 

\medskip
For the exceptional Lie algebras, Lemma \ref{lowb} shows that the Lie algebra
of type $G_2$ is generated by no fewer than 4 and the other four (types $F_4$
and $E_6$,  $E_7$,  $E_8$) are generated by no fewer than 5 extremal elements.
When $k$ has characteristic 3, the Lie algebra $\g$ of
Chevalley type $G_2$ has a 7-dimensional simple quotient.
This quotient
is isomorphic to a quotient of $\sli_3$
by its center, and is generated by 3 extremal elements.
Nevertheless, 
$\g$ itself, being 14-dimensional, cannot be generated by fewer
than four extremal elements.

\bigskip
To prove the theorem, it remains to show that there is a generating set of
extremal elements of the size indicated in the table.  To this end, we often
argue along the following pattern. Let $\g$ be the Lie algebra
under consideration.  We shall work with a fixed Chevalley basis $\Bc$ of $\g$
and a fixed root system whose simple roots $\alpha_1,\ldots,\alpha_n$ are
labeled as in \cite{Bou}. When we talk about root elements corresponding
to specific roots we mean elements from this basis.

We shall select a Lie subalgebra $M$ of $\g$ of Chevalley type
generated by extremal
elements, usually chosen from $\Bc$.  Let $\g = M \oplus V_1 \oplus \cdots
\oplus V_r$ be the decomposition of the $M$-module $\g$ into irreducible
$M$-modules.  The modules $V_i$ are often spanned by certain elements of $\Bc$
and usually the index $i$ corresponds to coefficients of simple roots not
supported in $M$.  We shall write down an element $d$ which is the image of an
extremal element under a composition of exponentials of root elements.
It is an extremal element whose projections onto
many of the modules $V_j$ $(j=1,\ldots,r)$ are nonzero. Next, we take $C$ to
be the Lie algebra generated by $M$ and $d$ and, by suitably bracketing $d$ by
elements of $C$, we find vectors from each $V_j$, usually again elements of
$\Bc$. Once those are found, we see that $C$ coincides with $\g$, so $t(\g)\le
t(M)+1$.  (In case $\g$ is of type $C_n$, we need two additional extremal
elements instead of 1 outside $M$ and there we show $t(\g) \le t(M)+2$.)     

We start with the classical Lie algebras. In this case, with exception of
$A_n$, we will write the simple roots using the unit orthonormal vectors in
$\CC^n$, $\eps_1, \dots, \eps_n$ with $\alpha_i$'s expressed in terms of
$\eps_i$'s as in Bourbaki.

All computations were made using the {\sc ELIAS} routines 
\cite{ELIAS} in GAP \cite{gap},
and {\sc LiE} \cite{Lie}.
The signs appearing here are the ones
arising from the implementations in the relevant packages.  In 
general, the signs depend on the choice of Chevalley basis.

\medskip

\subsection{Type $A_n$}

By Lemma \ref{ub} we know that $A_n$ can be generated by $n+1$ extremal
elements, that is, $t(A_n)\le n+1$. So there is nothing left to prove. 

\subsection{Type $B_n$}

Let $\g$ be the Lie algebra of type $B_n$. 

First assume $n=3$. Take $M$ to be the Lie subalgebra of $\g$ generated by the
root elements corresponding to the first two simple roots, $\eps_1 -
\eps_2$ and $\eps_2 - \eps_3$ and their negatives.
Then $M$ is of type $A_2$ and contains
the root elements $x_{\pm(\eps_1 - \eps_2)}$, $x_{\pm(\eps_2 - \eps_3)}$ and
$x_{\pm(\eps_1 - \eps_3)}$. Moreover, $\g$ decomposes as
\[ \g = M \oplus V_{-2} \oplus V_{-1}\oplus V_{1}\oplus V_{2}\oplus H_0,\]
where $V_1$ is the natural $M$-module linearly spanned by root elements
$x_{\eps_1}$,  $x_{\eps_2}$ and $x_{\eps_3}$; $V_2$ is spanned by root elements
$x_{\eps_1+\eps_2}$,  $x_{\eps_1+\eps_3}$ and $x_{\eps_2+\eps_3}$. Their dual
$M$-modules $V_{-1}$ and $V_{-2}$ are natural $M$-modules linearly spanned by
$x_{-\eps_1}$, $x_{-\eps_2}$, $x_{-\eps_3}$ and  $x_{-(\eps_1+\eps_2)}$,
$x_{-(\eps_1+\eps_3)}$, $x_{-(\eps_2+\eps_3)}$, respectively and $H_0$ is the
trivial $M$-submodule of $L$, spanned by the torus element $h_{\eps_3}$
centralizing $M$.

Now write:
\begin{multline*}
d =\exp(x_{-(\eps_1+\eps_2)}) \exp(x_{\eps_1})x_{-(\eps_1-\eps_2)} =\\ 
x_{-(\eps_1-\eps_2)} - x_{\eps_2}+x_{\eps_1+\eps_2} + x_{-\eps_1}+ h_{-(\eps_1+\eps_2)}- x_{-(\eps_1+\eps_2)}. 
\end{multline*}

Consider the Lie subalgebra $C$ of $\g$ generated by $M$ and $d$. As
$x_{\eps_1-\eps_2}\in M$, we have $- x_{\eps_1} - x_{-\eps_2}+M =
[x_{\eps_1-\eps_2}, d]+M \subseteq C$, whence $x_{\eps_1}+ x_{-\eps_2}\in
C$. Bracketing by the element $x_{\eps_2-\eps_3}$ of $M$, we find
$x_{-\eps_3}\in C$. In particular, the $M$-submodule $V_{-1}$ generated by
this element belongs to $C$, whence also $x_{-\eps_1}$, and so $e := d-
x_{-(\eps_1-\eps_2)} - x_{-\eps_1} \in C$. 

Taking brackets of $x_{-\eps_1}$ with $e$ gives $2 x_{-(\eps_1-\eps_2)}-
x_{\eps_2} +  x_{-\eps_1} \in C$, whence $x_{\eps_2} \in C$. In particular,
$x_{\eps_3}$, which lies in the same $M$-submodule, belongs to $C$, and hence
so do all root elements corresponding to positive roots. But then $e+C =
x_{-(\eps_1+\eps_2)} + C \subseteq C$, and so also the lowest root element
belongs to $C$. The root elements corresponding to the 
simple roots and to the
lowest root generate $\g$, so $C = \g$. Thus $t(\g)\le t(M)+1=4$. 

Next, assume $n\ge4$.  In the Lie algebra $\g$ of type $B_n$ we
have the Lie subalgebra $M$ of type $D_n$ generated by all long root elements
(extremal elements) of the form $x_{\eps_i\pm \eps_j}$ $(1\le i<j\le n)$.  As
an $M$-representation space, $\g$ decomposes into the adjoint module $M$ and
the natural representation $V_1$ of degree $2n$. Here, $V_1$ is linearly
spanned by the elements $x_{\pm\eps_j}$ $(j = 1,\ldots,n)$ from the Chevalley
basis $\Bc$. 

Now 
\[ d := \exp(x_{\eps_2})x_{\eps_1-\eps_2} = x_{\eps_1-\eps_2} - x_{\eps_1} -
x_{\eps_1+\eps_2}\] 
is an extremal element. Let $C$ be the Lie subalgebra of $\g$ generated by $M$
and $d$. 

We show it is $\g$. For, $d + M =  x_{\eps_1} + M$, so $x_{\eps_1}\in C$,
whence $V_1 \cap C \ne \{0\}$. By irreducibility of $V_1$, we find that
$x_{\pm\eps_j}\in C$ for all $j$. Thus, $C$ contains all root elements of the
standard Chevalley basis $\Bc$ of $\g$, and so coincides with $\g$. 

The conclusion is that $t(B_n) \le t(D_n)+1$. In particular, the proof of the
theorem for $B_n$ is complete once the theorem is shown to hold for $D_n$ (as
then $t(D_n)= n$). 

\subsection{Type $C_n$}

Suppose $n\ge2$ and let $\g$ be the Lie algebra of type
$C_n$. Take $M$ to be the Lie subalgebra of type $C_{n-1}$ generated by
the root elements with roots in the linear span of
$\eps_2,\ldots,\eps_n$. Then $\g$ decomposes as  
\[ \g = M\oplus V_1\oplus V_{-1}\oplus V_2\oplus V_{-2}\oplus V_{0},\]
where $V_1$ and $V_{-1}$ are natural $M$-modules, spanned by
all root elements with first simple
root coordinate $1$, respectively $-1$,
and $V_2$, $V_0$, $V_{-2}$ are trivial modules, spanned by $x_{2\eps_1}$,
$h_{2\eps_1}$, and $x_{-2\eps_1}$, respectively.  

Now take
\[ d_1 = \exp(x_{-(\eps_1+\eps_2)})x_{2\eps_1} = x_{2\eps_1} -
x_{\eps_1-\eps_2}-  x_{-2\eps_2}\] 
and
\[ d_2 = \exp(x_{\eps_1+\eps_2})x_{-2\eps_1} =
x_{-2\eps_1}+x_{-(\eps_1-\eps_2)} - x_{2\eps_2}.\] 

Let $C$ be the Lie subalgebra generated by $M$, $d_1$, and $d_2$. Then $d_1+M
= x_{2\eps_1} - x_{\eps_1-\eps_2} +M$ and $d_2 + M = x_{-2\eps_1} + x_{-(\eps_1-\eps_2)} +M$.  

Now $C = [x_{-2\eps_2},d_2] + C = x_{-(\eps_1+\eps_2)} + C$, and  $C =
[x_{2\eps_2},d_1] + C =  x_{\eps_1+\eps_2} + C$, so generators for each
natural $M$-submodule of $\g$ are in $C$. But then $C = d_1+C = x_{2\eps_1} +
C$ and  $C = d_2+C =x_{-2\eps_1} + C$, so $V_2$ and $V_{-2}$ are also
contained in  $C$. We have established that $\g\subseteq C$, whence $C = \g$.

The conclusion is $t(\g)\le t(M)+2$. If $n=2$, then $M$ is of type
$A_1$, so we have $t(M) = 2$ and $t(\g)\le 4$. By induction, we find $t(\g)\le
2n$ for arbitrary $n$.  

\subsection{Type $D_n$}\label{d4}

First we consider the case $n=4$. Let $M$ be the subalgebra of $\g$ generated
by the root elements corresponding to the roots $x_{\pm(\eps_1- \eps_2)}$,
$x_{\pm(\eps_2-\eps_3)}$ and $x_{\pm(\eps_1-\eps_3)}$. Then
$M$ is of type $A_2$ and $\g$ decomposes as
\[ \g = M \oplus V_{\eps_3-\eps_4}\oplus V_{\eps_3+\eps_4}\oplus V_{1,0}\oplus
V_{0,1}\oplus V_{1,1}\oplus V_{-1,0}\oplus V_{0,-1}\oplus V_{-1,-1},\]
where $V_{i,j}$ is a 3-dimensional module linearly spanned by the root
elements corresponding to simple roots of the form $\lambda
(\eps_1-\eps_2)+\mu (\eps_2-\eps_3) + i (\eps_3-\eps_4) +j (\eps_3+\eps_4)$
where $\lambda$, $\mu$ are scalars and $V_{\eps_3-\eps_4}$ and
$V_{\eps_3+\eps_4}$ are one-dimensional $M$-modules.

Consider
\begin{eqnarray*}
d &=& \exp(x_{\eps_1+\eps_4}) \exp(x_{-(\eps_3-\eps_4)})
\exp(x_{-(\eps_1+\eps_3)})x_{\eps_3-\eps_4} = \\
& & -x_{\eps_1+\eps_3} - 2 x_{\eps_1+\eps_4} + x_{\eps_3-\eps_4}- 2
x_{-(\eps_3-\eps_4)}\\
& & - x_{-(\eps_1+\eps_3)} + x_{-(\eps_1+\eps_4)} +
h_{-(\eps_3-\eps_4)}+ h_{\eps_1+\eps_4}. 
\end{eqnarray*}

This is an extremal element as it is the image of a root element under an
automorphism of $\g$.  

Denote by $C$ the Lie subalgebra of $\g$ generated by $M$ and $d$. A
computation shows 
\[ [x_{-(\eps_1-\eps_3)},[x_{\eps_1-\eps_2},d]] =  -
x_{-(\eps_2-\eps_3)}+x_{-(\eps_1+\eps_2)} .\]   
Since $x_{-(\eps_1-\eps_3)}$, $x_{\eps_1-\eps_2}$, $x_{-(\eps_2-\eps_3)}\in
M$, we derive $x_{-(\alpha_1+2\alpha_2+\alpha_3+\alpha_4)} =
x_{-(\eps_1+\eps_2)} \in C$; that is, the lowest root element belongs to
$C$. It follows that $V_{-1,-1}\subseteq C$. 

Another computation yields
\[ [x_{-(\eps_2-\eps_3)},[x_{-(\eps_1-\eps_2)},d]]  =  - 2 x_{\eps_3+\eps_4} +
x_{-(\eps_1-\eps_3)}.\] 
Similar to the above, we derive $x_{\alpha_4} = x_{\eps_3+\eps_4} \in C$. So
$V_{0,1}\subseteq C$.  


By computation, we have 
\[ [x_{\eps_1-\eps_2},[h_{\eps_1-\eps_3},[x_{-(\eps_1+\eps_3)},d]]] =
x_{-(\eps_2+\eps_4)}.\]  
Since $x_{\eps_1-\eps_2}$, $h_{\eps_1-\eps_3}\in M$ and
$x_{-(\eps_1+\eps_3)}\in V_{-1,-1}\subset C$, we see $x_{-(\alpha_2+\alpha_4)}
= x_{-(\eps_2+\eps_4)}\in C $, whence $V_{0,-1}\subset C$.

Since $x_{-(\eps_1+\eps_4)}\in V_{0,-1}\subseteq C$ and $h_{\eps_1-\eps_3}\in
M$, we have 
\[ C = [x_{-(\eps_1-\eps_3)},[x_{\eps_1-\eps_3},d]] +C = x_{-(\eps_1+\eps_4)}
+ 2 h_{\eps_1-\eps_3}+x_{\eps_3-\eps_4}+C =
x_{\eps_3-\eps_4}+C ,\]
and so $x_{\alpha_3} = x_{\eps_3-\eps_4}\in C$. So $V_{1,0}\subseteq C$.   

We have seen that the root elements corresponding to all four simple roots
and to the lowest root belong to $C$. Since they generate $\g$, we conclude
$C=\g $, whence $t(\g)\le t(M)+1 = 4$. 

\bigskip
Assume, from now on, that $n>4$. We show that the Lie algebra $\g$ of type
$D_n$ satisfies $t(\g) \le n$. 

Let $M$ be the Lie subalgebra of $\g$ generated by all root
elements corresponding to the subrootsystem on $\eps_i\pm \eps_j$ for $2\le
i<j\le n$. Then $M$ has type $D_{n-1}$.  As an $M$-module, $\g$ decomposes as
\[ \g = M\oplus V_1\oplus V_{-1}\oplus H_0,\]
where $H_0$ is a 1-dimensional toral subalgebra centralizing $M$ and $V_1$, $V_{-1}$ are
natural modules spanned by all root elements whose root vectors have first
coordinate (with respect to the basis of simple roots) equal to 1, $-1$,
respectively. 

By induction $t(M) \le n-1$, so we need to find an extremal element of $\g$
that, together with $M$, generates $\g$. Take 
\[ d = \exp(x_{-(\eps_1-\eps_2)})x_{\eps_1-\eps_2}= x_{\eps_1-\eps_2} -
h_{\eps_1-\eps_2}- x_{-(\eps_1-\eps_2)}. \]
Now $[x_{\eps_2-\eps_3},d] + M = -  x_{\eps_1-\eps_3} + M$, so
$x_{\eps_1-\eps_3}\in C$, whence every root element whose root has first root
coordinate 1 belongs to $C$. That is, $V_1 \subseteq C$ and therefore
$x_{\eps_1 -\eps_2} \in C$. Similarly, $C = [ x_{-(\eps_2-\eps_3)},d] + C =
- x_{-(\eps_1-\eps_3)}+C$, and so also every root element whose root has first
coordinate $-1$ belongs to $C$. That is, $V_{-1} \subseteq C$ and therefore
$x_{-(\eps_1 +\eps_2)} \in C$. Thus, all root elements from the standard Chevalley basis
belong to $C$, proving $C=\g$, and $t(\g)\le t(M)+1 = n$. 

\subsection{Type $E_6$}

Let us denote by $\alpha_1,\alpha_2,\alpha_3,\alpha_4,\alpha_5,\alpha_6$ the
simple roots of $\ek_6$. Let $M$ be a subalgebra of $\ek_6$ generated by $M =
\langle x_{\alpha_2}, x_{\alpha_3}, x_{\alpha_4}, x_{\alpha_5},
x_{-(\alpha_2+\alpha_3+  2\alpha_4+\alpha_5)} \rangle$. Then $M$ is of type
$D_4$.   
One has the following decomposition of $\ek_6$:
\[ \ek_6 = M \oplus W \oplus V_{0,1} \oplus V_{1,0} \oplus
V_{1,1} \oplus V_{0,-1} \oplus V_{-1,0} \oplus V_{-1,-1}, \]
where $V_{a, b}$ is the module generated by all $x_{\alpha}$ with $a$
coefficient of $\alpha_1$ and $b$ coefficient of $\alpha_6$ in $\alpha$. The
action of $M$ on $W$ is trivial.

We need one more extremal element:
\begin{eqnarray*}
d&=&
\exp(x_{-(\alpha_1+\alpha_3+\alpha_4)})\exp(x_{\alpha_3+\alpha_4+\alpha_5+\alpha_6})\exp(x_{-(\alpha_1+\alpha_2+\alpha_3+2\alpha_4+2\alpha_5+\alpha_6)})
x_{\alpha_1} \\
&=& x_{\alpha_1} - x_{-(\alpha_2+\alpha_3+2\alpha_4+2\alpha_5+\alpha_6)} + 
x_{\alpha_1+\alpha_3+\alpha_4+\alpha_5+\alpha_6} +
x_{-(\alpha_2+\alpha_4+\alpha_5)}\\
& & -  x_{-(\alpha_3+\alpha_4)} 
- x_{-(\alpha_1+\alpha_2+2\alpha_3+3\alpha_4+2\alpha_5+\alpha_6)}\\
& &+
x_{\alpha_5+\alpha_6} - x_{-(\alpha_1+\alpha_2+\alpha_3+2\alpha_4+\alpha_5)}.
\end{eqnarray*}

We claim that the subalgebra $C$ generated by $M$ and $d$ is the whole algebra
$\ek_6$. 

The root elements $x_{-(\alpha_2+\alpha_4+\alpha_5)}$ and
  $x_{-(\alpha_3+\alpha_4)}$ belong to $M$, so $e:= d -
  x_{-(\alpha_2+\alpha_4+\alpha_5)} +  x_{-(\alpha_3+\alpha_4)}$ is in $C$.    
Since $x_{\alpha_3}$, $x_{-\alpha_5}$ and $x_{- \alpha_2}$ are in $M$,
the following brackets are in $C$:
\[ [x_{\alpha_3}, e] =  x_{\alpha_1+\alpha_3}, \; \; \; [x_{-\alpha_5}, e] =
 x_{\alpha_6}  \text{  and  } [x_{-\alpha_2}, e] = -
 x_{-(\alpha_1+2\alpha_2+2\alpha_3+3\alpha_4+2\alpha_5+\alpha_6)} \] 
But $x_{\alpha_1+\alpha_3} \in C$ implies $V_{1,0} \subset C$, and so
$x_{\alpha_1} \in C$. 

Now $C$ contains all the simple root elements and the lowest root
element. By Lemma \ref{gen}, we are done.

\subsection{Type $E_7$}

Let us denote by
$\alpha_1,\alpha_2,\alpha_3,\alpha_4,\alpha_5,\alpha_6,\alpha_7$ the
simple roots of $\ek_7$. Take the subalgebra $M$ of type $D_4$ in
$\ek_7$ exactly as in $\ek_6$.

To generate the whole algebra $\ek_7$, consider the extremal element
\begin{eqnarray*}
d&=&
\exp(x_{\alpha_2+\alpha_3+\alpha_4+\alpha_5+\alpha_6+\alpha_7})
\exp(x_{-(\alpha_1+\alpha_3)})\\
& &\exp(x_{-(\alpha_1+\alpha_2+\alpha_3+2\alpha_4+\alpha_5+\alpha_6)})
\exp(x_{\alpha_3+\alpha_4+\alpha_5+\alpha_6})\\
& &\exp(x_{-(\alpha_1+\alpha_2+\alpha_3+2\alpha_4+2\alpha_5+\alpha_6+\alpha_7)})x_{\alpha_1}\\
&=& x_{\alpha_1} +  x_{\alpha_4+\alpha_5+\alpha_6} +
x_{\alpha_1+\alpha_3+\alpha_4+\alpha_5+\alpha_6} + 
x_{\alpha_2+\alpha_4+\alpha_5+\alpha_6+\alpha_7}\\
& &+ x_{\alpha_1+\alpha_2+\alpha_3+ \alpha_4+\alpha_5+\alpha_6+\alpha_7}-
x_{\alpha_1+\alpha_2+2\alpha_3+ 3\alpha_4+2\alpha_5+\alpha_6}\\
& & + x_{\alpha_1+2\alpha_2+2\alpha_3+3 \alpha_4+2\alpha_5+\alpha_6+\alpha_7} - 
x_{-\alpha_3}+ x_{-(\alpha_2+\alpha_4)} - x_{-(\alpha_4+\alpha_5)}\\
& & -
x_{-(\alpha_1 +\alpha_2+ \alpha_3+\alpha_4)} + x_{-(\alpha_1 +
  \alpha_3+\alpha_4+\alpha_5)} - 
x_{-(\alpha_2 +\alpha_3+2\alpha_4+\alpha_5+\alpha_6)} \\
& &-x_{-(\alpha_1+\alpha_2
  +2\alpha_3+2\alpha_4+\alpha_5+\alpha_6)} - x_{-(\alpha_2+\alpha_3
  +2\alpha_4+2\alpha_5+\alpha_6+\alpha_7)} -\\
& & x_{-(\alpha_1+\alpha_2+2\alpha_3+2 \alpha_4+2\alpha_5+\alpha_6+\alpha_7)}  -
x_{-(\alpha_1+2\alpha_2+2\alpha_3 +4
  \alpha_4+3\alpha_5+2\alpha_6+\alpha_7)}+\\ 
& &x_{-(2\alpha_1+2\alpha_2+3\alpha_3 +4\alpha_4+3\alpha_5+2\alpha_6+\alpha_7)}.
\end{eqnarray*}
We show that the subalgebra $C$ generated by $M$ and $d$ is the whole
algebra  $\ek_7$. 
The root elements $x_{-\alpha_3}$, $x_{-(\alpha_2+\alpha_4)}$ and
$x_{-(\alpha_4\alpha_5)}$ belong to $M$, so $e := d + x_{-\alpha_3} -
x_{-(\alpha_2+\alpha_4)} +  x_{-(\alpha_4+\alpha_5)} $ is in $C$.    
The Lie algebra $\ek_7$ decomposes as:
\begin{eqnarray*}
\ek_7 &=& M \oplus W \oplus V_{0,0,1} \oplus V_{0,1,0}
\oplus V_{0,1,1} \oplus V_{0,2,1} \oplus V_{1,0,0} \oplus V_{1,1,0} \\
& &\oplus
V_{1,1,1} \oplus V_{1,2,1} \oplus V_{2,2,1} \oplus V_{0,0,-1} \\
& & \oplus V_{0,-1,0} \oplus V_{0,-1,-1} \oplus V_{0,-2,-1} \\
& & \oplus V_{-1,0,0} \oplus V_{-1,-1,0} \oplus V_{-1,-1,-1} \\
& &\oplus V_{-1,-2,-1} \oplus V_{-2,-2,-1},
\end{eqnarray*}
where $W$ is a 3-dimensional module with trivial $M$ action. The modules
$V_{1,0,0}$, $V_{0,1,1}$, $V_{0,1,0}$ and $V_{1,2,1}$ have dimension 8 and
$V_{a,b,c}$ is generated by all roots elements not in M such that $a$ is the
coefficient of $\alpha_1$, $b$ of $\alpha_6$ and $c$ of $\alpha_7$. 

The bracket
$[x_{\alpha_2+\alpha_3+2\alpha_4+\alpha_5},[x_{\alpha_3}, e]]$ is easily computed
to be a nonzero scalar multiple of
$x_{\alpha_1+\alpha_2+2\alpha_3+2\alpha_4+\alpha_5}$, and,
since both arguments of the bracket belong to $M$, the element
$x_{\alpha_1+\alpha_2+2\alpha_3+2\alpha_4+\alpha_5}$ belongs to $C$.
Consequently, $V_{1,0,0}$
is contained in $C$ and so $x_{\alpha_1} \in C$. 

Since $x_{-(\alpha_4+\alpha_5)}, x_{\alpha_5} \in C$, we have
$[x_{\alpha_5},[x_{-(\alpha_4+\alpha_5)}, e]] = - x_{\alpha_5+\alpha_6} \in C$.
So $V_{0,1,0} \subset C$ and $x_{\alpha_6} \in C$.
Also, $[x_{-\alpha_5},[x_{\alpha_3}, e]]$ belongs to $C$. Expanding it,
we find
$x_{-(\alpha_1+\alpha_2+\alpha_3+2\alpha_4+2\alpha_5+\alpha_6)} \in C$.  Hence
$C$ contains $V_{-1,-1,0}$ and taking brackets with $x_{\alpha_1}$ we obtain 
$V_{0,-1,0} \subset C$. 

Now $x_{\alpha_7} \in C$, as
$[x_{-(\alpha_3+\alpha_4+\alpha_5+\alpha_6)},[x_{-(\alpha_1+\alpha_2+\alpha_3+\alpha_4+\alpha_5+\alpha_6)},e]]$ is a nonzero element of 
$k x_{\alpha_7} +  V_{0,-1,0}$ which also belongs to $C$.


Finally, as $x_{\alpha_5+\alpha_6}\in C$ and
$ [x_{-\alpha_4}, e] =  x_{\alpha_5+\alpha_6} -
x_{-(\alpha_1+\alpha_2+2\alpha_3+3\alpha_4+2\alpha_5+\alpha_6+\alpha_7)}$,
we find  $V_{-1,-1,-1} \subset C$.
Taking brackets with a convenient element of $V_{-1,-1,0}$, we see that also the
lowest root element is in $C$. Now $C=\ek_7$ by Lemma \ref{8.1}.

\subsection{Type $E_8$}

Let $\alpha_1,\alpha_2,\alpha_3,\alpha_4,\alpha_5,\alpha_6,\alpha_7,\alpha_8$
be the simple roots of $\ek_8$. Adopt here the notation $(a_1, a_2, a_3, a_4,
a_5, a_6, a_7, a_8)$ for a root $\alpha = \sum_{i=1}^8 a_i \alpha_i$.
Consider the subalgebra $M$ of type $D_4$ in $\ek_8$ as in $\ek_6$.  
We have the following decomposition of $\ek_8$ into $M$-modules.
$$\ek_8 = M \oplus W \oplus
 \bigoplus_{a,k,l,m \in I} V_{a,k,l,m} ,$$
where $W$ is a 28 dimensional trivial module,
$V_{a,k,l,m}$ is an irreducible 8 dimensional module generated by all roots
elements $x_{\alpha}$ where $a$, $k$, $l$, $m$ are respectively the
coefficients of $\alpha_1$, $\alpha_6$, $\alpha_7$, $\alpha_8$ in $\alpha$,
and
\begin{eqnarray*}
I&=& \pm \{(0,1,0,0), (0,1,1,0), (0,1,1,1), (1,0,0,0), (1,1,0,0),\\
&&\quad\  (1,1,1,0), (1,1,1,1), (1,2,1,0), (1,2,1,1), (1,2,2,1),\\
&&\quad\  (1,3,2,1), (2,3,2,1)\}. 
\end{eqnarray*}
To generate the whole algebra $\ek_8$, consider the extremal element
\begin{eqnarray*}
d&=&\exp(x_{01111110})\exp(x_{-11121110})\exp(x_{01122111})\exp(x_{-10111111})\\
& &\exp(x_{12343321})\exp(x_{-23354321})x_{10000000}\\
&=&
x_{10000000}- x_{01011000}+ x_{01111100} \\
& &+x_{00111111} +x_{11111110} + x_{11122111}- x_{11222211}-
x_{12232210}\\
& &-x_{12233221}- x_{22343321}-x_{13354321}- x_{23465432}+
x_{-10000000}  \\
& &-x_{-00010000} - x_{-01011000} - x_{-00111111} + x_{-11121100}
-x_{-01121110}\\
& &-x_{-11122111}- x_{-12232210}-x_{-11232221}+
x_{-12243211}+x_{-22343321}\\
& &+x_{-13354321}- x_{-23465432}.
\end{eqnarray*}

Let $C$ be the Lie subalgebra generated by $M$ and $d$. 
Bracketing and using the decomposition of $\ek_8$ as an $M$-module,
one can again derive
 equality between $C$ and $\ek_8$. We omit further details,
as they are very similar to the other cases dealt with (but lengthier).

\subsection{Type $F_4$}

Write $\alpha_1$, $\alpha_2$, $\alpha_3$ and $\alpha_4$ for the simple
roots of
$\fk_4$.  So, $\alpha_1 = \eps_2 - \eps_3$, $\alpha_2 =
\eps_3 - \eps_4$, $\alpha_3 = \eps_4$, $\alpha_4 = \frac{1}{2}
\left(  \eps_1 - \eps_2- \eps_3 - \eps_4 \right)$.

Let $M$ be the subalgebra of $\fk_4$ generated by $x_{\alpha_1}$,
  $x_{\alpha_2}$, $x_{\alpha_2 + 2 \alpha_3}$, $x_{\alpha_2 + 2 \alpha_3+ 2
  \alpha_4}$ and $x_{-(2\alpha_1+3\alpha_2 + 4 \alpha_3+ 2 \alpha_4)}$. That
  is, $M = \langle x_{\eps_1 - \eps_2},x_{\eps_2 - \eps_3},x_{\eps_3
  - \eps_4},x_{\eps_3  + \eps_4}, x_{ -(\eps_1+ \eps_2)}\rangle$. According to
  Lemma \ref{gen}, $M$ is of type $D_4$. By Section \ref{d4} we know then that
  $M$ is generated by four extremal elements.  

Define the following element in $\fk_4$:
\begin{multline*}
d =
\exp(x_{-(\alpha_1+2\alpha_2+3\alpha_3+\alpha_4)})\exp(x_{\alpha_4})(x_{\alpha_2+2\alpha_3})
=  x_{\alpha_2+2\alpha_3} - x_{\alpha_2+2\alpha_3+\alpha_4}-\\
x_{\alpha_2 + 2\alpha_3+2\alpha_4} - x_{-(\alpha_1 +  \alpha_2+\alpha_3)} -
x_{-(\alpha_1+\alpha_2+\alpha_3+\alpha_4)} - x_{-(2\alpha_1+3\alpha_2 + 4  \alpha_3+ 2 \alpha_4)}.
\end{multline*}
The element $d$ is extremal since $\alpha_2+2\alpha_3$ is long. Consider the
subalgebra $C$ generated by $M$ and $d$. 

The root elements $x_{\alpha_2+2\alpha_3}$, $x_{\alpha_2 + 2 \alpha_3+ 2
  \alpha_4}$  and $x_{-(2\alpha_1+3\alpha_2 + 4 \alpha_3+ 2 \alpha_4)}$ are in
  $M$. So $f := d - x_{\alpha_2+2\alpha_3}+ x_{\alpha_2 + 2 \alpha_3+ 2
  \alpha_4}  +  x_{-(2\alpha_1+3\alpha_2 + 4 \alpha_3+ 2 \alpha_4)} =
  -x_{\alpha_2 + 2 \alpha_3+ \alpha_4} - x_{-(\alpha_1 +  \alpha_2+\alpha_3)}
  - x_{-(\alpha_1+\alpha_2+\alpha_3+\alpha_4)}$ is in $C$. 

The Lie algebra $\fk_4$ decomposes as:
\begin{eqnarray*}
\fk_4 &=& M \oplus V_{0,1} \oplus V_{1,0} \oplus V_{0,0} \oplus V_{1,1} \oplus
V_{1,2}\\
& & \oplus V_{0,-1} \oplus V_{-1,0} \oplus W_{0,0} \oplus V_{-1,-1} \oplus
V_{-1,-2},
\end{eqnarray*}
where $V_{0,1}$ is the module generated by $x_{\alpha_4}$, $x_{\alpha_3 +
  \alpha_4}$, $x_{\alpha_2 + \alpha_3+\alpha_4}$ and $x_{\alpha_2 +
  2\alpha_3+\alpha_4}$, $V_{1,0} =\langle
x_{\alpha_1+\alpha_2+\alpha_3}\rangle$, $V_{0,0} =\langle x_{\alpha_3},
x_{\alpha_2 + \alpha_3} \rangle$, $V_{1,1} =\langle
x_{\alpha_1+\alpha_2+\alpha_3+\alpha_4}, x_{\alpha_1+\alpha_2 +2
  \alpha_3+\alpha_4}, \linebreak x_{\alpha_1+2\alpha_2+2\alpha_3+\alpha_4},
x_{\alpha_1+2\alpha_2 +3 \alpha_3+\alpha_4} \rangle$ and $V_{1,2} =\langle
x_{\alpha_1+2\alpha_2+3\alpha_3+2\alpha_4} \rangle$. The module $W_{0,0}$ is
generated by $x_{-\alpha_3}$, $x_{-(\alpha_2 + \alpha_3)}$ and obviously the
remaining modules are generated by the negative root elements. 

We have $[x_{\alpha_1 +\alpha_2 + 2 \alpha_3}, f] + C =  x_{\alpha_3}  + C$
since $x_{\alpha_1 +\alpha_2 + 2  \alpha_3} \in M$. Hence $x_{\alpha_3}$ is in
$C$ and then $V_{0,0} \subset C$. 

Now the bracket of $x_{\alpha_1 +\alpha_2 + 2\alpha_3+ 2 \alpha_4} \in M$ with
$f$ leads to the conclusion that $x_{\alpha_3+\alpha_4} \in C$, thus $V_{0,1}
\subset C$. In particular, $x_{\alpha_4} \in C$.

By Lemma \ref{gen}, the algebra $\fk_4$ is generated by $x_{\alpha_1}$,
  $x_{\alpha_2}$, $x_{\alpha_3}$, $x_{\alpha_4}$ and \linebreak
  $x_{-(2\alpha_1+3\alpha_2   + 4 \alpha_3+ 2 \alpha_4)}$. But $C$ contains
  all these root   elements. Therefore $C = \fk_4$. The conclusion is
  $t(\fk_4)\le t(D_4)+1 = 5$. 

\subsection{Type $G_2$}

Let $\alpha$ and $\beta$ be the simple roots of $\g_2$. Let $M$ be the Lie
subalgebra of $\g_2$ generated by the long root elements $x_{\beta}$,
$x_{3\alpha+\beta}$,  $x_{-(3\alpha+2\beta)}$. This subalgebra is of type
$A_2$. 

Consider the extremal element $d = \exp(x_{-(2\alpha+\beta)})
(x_{3\alpha+2\beta})$; it is extremal since $3\alpha+2\beta$ is a long
root. We have
\[ d = x_{3\alpha+2\beta} + x_{\alpha+\beta} -  x_{-\alpha}. \]
Let $C$ be the subalgebra of $\g_2$ generated by $M$ and $d$. We will prove that $C$
is actually $\g_2$. Since $x_\alpha$, $x_{\beta}$ and $x_{-(3\alpha+2\beta)}$
generate $\g_2$ and the latter two are obviously in $C$, it is enough to show
that $x_{\alpha}$ is also in $C$. 

Since $d - x_{3\alpha+2\beta} \in C$, so does the bracket 
\[  [x_{-(3\alpha +2\beta)}, d - x_{3\alpha+2\beta}] = [x_{-(3\alpha+2\beta)},
 x_{\alpha+\beta} - x_{-\alpha}] = - x_{-(2\alpha+\beta)}. \]

Therefore $[x_{3\alpha+\beta}, x_{-(2\alpha+\beta)}] = x_{\alpha}$ is also in
$C$.

As before we conclude that $C = \g_2$ and so $\g_2$ is generated by four
extremal elements.

\section{The bilinear form and the Killing form}
\label{bilsec}

Throughout this section, $L$ is a finite dimensional Lie algebra
(over the field $k$ of characteristic not $2$) generated by extremal elements.
Recall that the Killing form
$\kappa$ is defined by $\kappa(x,y) := \tr(\ad_x\ad_y)$ for $x,y \in L$.
We consider connections between the associative bilinear form $f$
defined in Theorem \ref{3.6} 
and the Killing form $\kappa$.

Fix $x \in \Ec$ and $y \in L$. Set $\phi := \ad_x \ad_y$.

\begin{Lm}  \label{cal}
We have $\phi^2 + \frac{1}{2}f(x,y) \phi: L \to kx + k[x,y]$.
\end{Lm}

\begin{proof}
For $z \in L$, using Lemma \ref{ide}, we find
$$\phi^2(z) = [x,[y,[x,[y,z]]]] = 
\frac{1}{2}(f(x,[y,[y,z]])x - f(x,[y,z])[x,y] - f(x,y)[x,[y,z]]).$$
This implies that $\phi^2+\frac{1}{2} f(x,y) \phi$ maps $L$ to  $kx+k[x,y]$.
\end{proof}

\begin{Lm}\label{lem}\label{9.2}
Let $x \in \Ec$, $ y \in L$. Then $\phi$ satisfies the following properties.
\begin{enumerate}
\item[a)] If $f(x,y) = 0$, then all eigenvalues of $\phi$ are $0$.
In particular, $\kappa(x,y) = 0$.
\item[b)] If $f(x,y) = -2$, then $\phi$ has eigenvalue $2$ with multiplicity
$2$, eigenvalue $1$ with multiplicity $s - 2$ and the remaining
eigenvalues are $0$. Here $s = \dim \ad_{x}(L)$.
In particular, $\kappa(x,y) = s + 2$.
\end{enumerate}

%
\end{Lm}

\begin{proof}
Note that $\phi(x)=-f(x,y)x$ and (by Lemma~\ref{ide})
$\phi([x,y])=-f(x,y)[x,y]$.
Hence $U = kx + k[x,y]$ is invariant under $\phi$.
Furthermore, $\phi^2 + \frac{1}{2} f(x,y) \phi$ is zero on $L/U$ by
Lemma~\ref{cal}. We write $\phi$ on $L/U$ in Jordan normal form with
eigenvalues $0$ or $-\frac{1}{2} f(x,y)$. Now a) is obvious.
Next suppose $f(x,y) = -2$. Then $x, [x,y]$ are linearly independent.
Denote by $E$
the eigenspace of $\phi$ (with eigenvalue 1).
Then $E+U/U = \ad_x(L)/U$ and $\dim E = \dim \ad_x(L) - 2$.
This proves the lemma.
\end{proof}

\bigskip
We now consider connections between the radicals of the two forms.
Recall that $\Rad(f)$ denotes the radical of $f$;
similarly, we write $\Rad(\kappa)$ for the radical of $\kappa$.

\begin{Cor}\label{rad}
We have $\Rad(f) \subseteq \Rad(\kappa)$.
\end{Cor}

\begin{proof}
Let $y \in \Rad(f)$. Then $f(x,y) = 0$ for all $x \in \Ec$. Hence
$\kappa(x,y) = 0$ for all $x \in \Ec$ by Lemma~$9.2$.  Therefore $y \in
\Rad(\kappa)$.   
\end{proof}

\begin{Remark}
Lemma \ref{lem} gives a way to tell if one of the
Lie algebras of Chevalley type is the radical of its
Killing form.
(In this case $\Rad(f)$ is properly contained in $\Rad(\kappa)$.)
For a root $\alpha$, set $\gamma_\alpha=
\dim \ad_{x_\alpha}( L) - 2$. Then by Lemma \ref{9.2},
$\kappa(x_\alpha,x_{-\alpha})= \gamma_\alpha+4$.
When there is only one root length,
$\gamma_\alpha$ is the number of roots which have inner product $-1$
with $\alpha$.
We look at the case of $\sli_3$. For $\alpha = e_i - e_j$,
there are exactly two such roots, $e_k - e_i$, and $e_j - e_k$,
so $\kappa(x_\alpha,x_{-\alpha})= 2+4 = 6$.
Hence $\sli_3$ has Killing form zero if $k$ is of characteristic $3$.
\end{Remark}

\medskip
\begin{Lm}
If $\ch(k) = 0$, then $\Rad(f) = \Rad(\kappa)$.
\end{Lm}

\begin{proof}
One inclusion comes from Corollary \ref{rad}. Let $y \in \Rad(\kappa)$. 
Assume that $y \not \in \Rad(f)$. Then there exists $x \in \Ec$ such that
$f(x,y) = -2$. Hence $\kappa(x,y) = \dim \ad_x(L) + 2 \neq 0$
in characteristic $0$.  This is a contradiction.
\end{proof}

\bigskip
In what follows we determine some more properties about the radicals of the
Killing form and of $f$.

\begin{Lm}\label{aa}
Let $J$ be an ideal of $L$ and $N := \spa_{k} \{ x \in \Ec \mid x \notin J
\}$. Then  $f(N,J) = 0$.
\end{Lm}

\begin{proof}
It is enough to show for $x \in \Ec \setminus J$ and $y \in J$ that 
$f(x,y)=0$.
We have $f(x,y) x = [x,[x,y]] \in J $. Hence $f(x,y) = 0$.
\end{proof} 

\medskip
\begin{Lm}\label{bb}
Suppose that $K$ is a solvable ideal of $L$.
Then $\Ec \cap K \subseteq \Rad(f)$.
\end{Lm}

\begin{proof}
Let $x \in \Ec \cap K$. Take $y \in \Ec$. If $y \notin K$, then $f(x,y) = 0$
by  Lemma \ref{aa}. If $y \in K$, suppose $f(x,y) \neq 0$. Then
$\langle x, y \rangle \simeq \sli_2$. Hence $\langle x,y \rangle$ is a 
non-solvable Lie  subalgebra of $K$ which is a contradiction. So $f(x,y) = 0$.

This means that $f(x,y) = 0$ for all $x \in \Ec \cap K$ and $y \in \Ec$.
The lemma follows as $\Ec$ linearly spans $L$, cf.\ Lemma \ref{linspan}.
\end{proof}

\bigskip
Let $\Rad(L)$ be the maximal solvable ideal in $L$. 

\begin{Prop}\label{ppp}
We have $\Rad(L)\subseteq \Rad(f)$.
\end{Prop}

\begin{proof}
Suppose that $K$ is a solvable ideal of $L$ 
and let $x \in K$. For $y \in \Ec$, either $y \in \Ec \setminus K$, 
in which case $f(x,y) =
0$ by Lemma \ref{aa}, or $y \in \Ec \cap K $ and then $f(x,y) = 0$ by Lemma \ref{bb}.
\end{proof}

\begin{Remark}\label{9last}
An example where $\Rad(L)$ is properly contained in $\Rad(f)$ is the
14-dimensional Lie algebra $L$ of type $G_2$ over a field $k$ of
characteristic 3.  Then, $\Rad(L) = 0$ but $\Rad(f)$ is a
7-dimensional ideal, generated by the short root elements of a
Chevalley basis.  Proposition \ref{radfsol} shows that characteristic
3 is indeed exceptional.

In fact, $\Rad(f)$ and $\Rad(L)$ are just two ideals of a chain of five.
Denote by $\SanRad(L)$
the linear span of all {\em sandwiches} of $\Ec$,
that is, those $x\in \Ec$, for which $\ad_x^2 = 0$.
We claim that $\SanRad(L)$ is an ideal of $L$.
Indeed,
if $x\in \SanRad(L)\cap \Ec$ and $y\in \Ec$, then $f_x=0$, so Corollary \ref{3.8}
shows that either $[x,y]=0$ or $[x,y]\in\Ec$ with $f_{[x,y]}=0$,
that is, $[x,y]\in \SanRad(L)$.
Since the restriction of $f$ to $\SanRad(L)$ vanishes,
the latter is a nilpotent Lie subalgebra of $L$ by Lemma \ref{4.2}, so
$\SanRad (L)$ is contained in $\NilRad(L)$,
the nilpotent radical of $L$.
We thus have the chain
$$\SanRad(L)\subseteq \NilRad(L)\subseteq
 \Rad(L)\subseteq \Rad(f)\subseteq \Rad(\kappa).$$
We know that the last two inclusions may be proper, and, in view Case (2)
of Theorem \ref{3gen},
we have examples where $\SanRad(L)$ is strictly contained in $\Rad(L)$.
\end{Remark}

\begin{Lm}\label{gabor}
Suppose $x\in\Ec \setminus\Rad(f)$ and $y\in \Rad(f)$.
Then $[x,y]$ satisfies $\ad_{[x,y]}^4 = 0$.
\end{Lm}

\begin{proof}
Put $X=\ad_x$ and $Y=\ad_y$. Then,
by \eqref{ide2} of Lemma \ref{ide}, for $z\in L$,
\begin{eqnarray*}
2XYXYz &= & 2[x,[y,[x,[y,z]]] \\
&=& f(x,[y,[y,z]])x - f(x,[y,z])[x,y] - f(x,y)[x,[y,z]] = 0.
\end{eqnarray*}
For the last equality, use that $f$ is associative and $y\in\Rad(f)$.
We obtain $XYXY=0$.
Also $X^2Y=0$, as $X^2Y(L)\subseteq \Rad(f) \cap kx = 0$.
Consequently,
\begin{eqnarray*}
\ad_{[x,y]}^4 &=& (XY-YX)^4 = 
(XYXY-XY^2X-YX^2Y+YXYX)^2  \\
&=& (YXYX-XY^2X)^2  \\
&=& XY^2X^2Y^2X+(YX)^4-XY^2(XY)^2X-YXYX^2Y^2X 
=  0.
\end{eqnarray*}
\end{proof}

\begin{Lm}  \label{ovf}
Let $K$ be an ideal of $L$ which is contained in $\Rad(f)$.
Then $\ov{L} := L/K$ is linearly spanned by extremal elements
with induced form $\ov{f}$ defined by $\ov{f}(\ov{x}, \ov{y})
:= f(x,y)$ for $x , y \in L$.
\end{Lm}

\begin{proof}
Since $K \subseteq \Rad(f)$, the expression 
$\ov{f}(\ov{x}, \ov{y}) := f(x,y)$ is well-defined for $x , y \in L$.
For $x \in \Ec$, $y \in L$, $[\ov{x}, [\ov{x}, \ov{y}]] =
\ov{f}(\ov{x}, \ov{y})\cdot \ov{x}$, whence $\ov{x} \in \Ec(\ov{L})$.
%
\end{proof}

\bigskip
We owe the proof of the following result to Gabor
Ivanyos.

\begin{Prop}\label{radfsol}
If the characteristic of $k$ is not 2 or 3, then $\Rad(f) = \Rad(L) $.
\end{Prop}

\begin{proof}
Recall that $\Rad(L)\subseteq \Rad (f)$, so if $\Rad(f)$ is contained in the center $ Z(L)$ of $L$,
there is nothing to prove.
Suppose, therefore, $\Rad(f)\not\subseteq Z(L)$.
We show that $\SanRad(L)\ne0$.

By Lemma \ref{gabor}, there is a nonzero element $y\in \Rad(f)$ with
$\ad_y^4 = 0$.  (If $\Ec\cap \Rad(f)\ne \emptyset$, any element of the
intersection will do, otherwise, take $y\in\Rad(f)$, and $x\in\Ec$
such that $[x,y]\ne0$, whose existence is guaranteed by the hypotheses
that $L=\langle\Ec\rangle$ and $\Rad(f)\not\subseteq Z(L)$, and apply
Lemma \ref{gabor}.)  By Proposition 2.1.5 of \cite{Kostr} (see also
Proposition 1.5 of \cite{Ben}), if the
characteristic of $k$ is not 2 or 3, the element $z = \ad_y^3(x)$ for
any $x\in L$, satisfies $\ad_z^3 = 0$.  In particular, there is a
nonzero element $y\in \Rad(f)$ with $\ad_y^3=0$.

%

If $\ad_y^2=0$, then we are done. Otherwise, there is
$b\in \Ec$ with
$x = \ad_y^2(b)\ne0$.

If $b\in \Rad(f)$, then $b\in \Rad(f)\cap\Ec\subseteq \SanRad(L)$, 
and we are done. So,  assume $b\not\in \Rad(f)$.
By Lemma 1.7(iii) of \cite{Ben}, as $k$ has characteristic not $2$ or
$3$, we have
$\ad_x^2 = \ad_y^2 \ad_b^2\ad_y^2$.
Since  $y\in \Rad(f)$, we have $\ad_y^2(L)\subseteq \Rad(f)$ and since
$b\in \Ec\setminus\Rad(f)$, we have $\ad_b^2(\Rad(f))=0$,
whence $\ad_x^2 = 0$, proving that $\Ec\cap \Rad(f) \ne 0$.
This establishes $\SanRad(L)\ne0$.

Thus, if $\Rad(f) \not\subseteq Z(L)$, it contains a nonzero sandwich,
and so $\NilRad(L)\ne0$. But then, in view of Lemma \ref{ovf}
and by induction on the dimension,
$\Rad(f)$ is solvable.

\end{proof}

\bigskip
We finish this section with a result identifying $\Rad(f)$ and
$\Rad(L)$ in arbitrary characteristic distinct from 2 under additional
hypotheses on the structure of $L$.
It shows that $f$ plays a role similar to $\kappa$ in the theory of Lie algebras
of characteristic $0$.

\begin{Lm} \label{orthdecomp}
Let $L =  L_1 \oplus L_2$ be a direct sum of ideals.
Then $L_1$ and $L_2$ are linearly spanned by extremal
elements (with form $f$ restricted to $L_1$ and $L_2$, respectively).
Furthermore, $L_1$ and $L_2$ are orthogonal with respect to $f$.
\end{Lm}

\begin{proof}
As $[L_1,L_2] = 0$, 
we have $\Ec(L_i) \subseteq \Ec$ for $i = 1,2$.
Let $x \in \Ec(L)$ and $y \in L$.
Write $x = x_1 + x_2$, $y = y_1 + y_2$ with $x_1, y_1 \in L_1$,
$x_2, y_2 \in L_2$.
We calculate $f(x,y)x_1 + f(x,y)x_2 = [x_1, [x_1,y]] + [x_2, [x_2, y]]$.
Hence $[x_i,[x_i,y]] = f(x,y)x_i$ for all $y \in L$, and
$x_i \in \Ec$ $(i=1,2)$.

Since $L$ is linearly spanned by $\Ec$, each $L_i$ is linearly
spanned by the projections. Hence $L_i$ is linearly spanned by
extremal elements (with form $f$ restricted to $L_i$).

Finally, for $z \in \Ec(L_1)$, $l_2 \in L_2$, we have
$f(z,l_2) z = [z, [z,l_2]] = 0$.
Since $L_1$ is linearly spanned by extremal elements, this shows that
the decomposition $L = L_1 \oplus L_2$ is an orthogonal one with respect
to $f$.
\end{proof}

\begin{Prop}\label{9.14}
We have
 $\Rad(f) = 0$
if and only if
$L$ is
a direct sum of simple ideals. 
\end{Prop}

\begin{proof}
Assume that $L = L_1 \oplus L_2 \oplus \dots \oplus L_r$ is 
a direct sum of simple ideals $L_i$.
By Lemma~\ref{orthdecomp}, the decomposition is an orthogonal one with
respect to $f$.
Suppose $L_i \subseteq \Rad(f)$ for some $i$. 
Then the form $f$ restricted to $L_i$ is trivial.
But then $L_i$ is nilpotent by Lemma~\ref{finite1}, 
a contradiction.
This means that the form $f$ restricted to each $L_i$ has a trivial
radical. Hence $\Rad(f) = 0$. 

As for the converse, assume $\Rad(f) = 0$.  Suppose that $I$ is a minimal
nonzero ideal of $L$. By Proposition \ref{ppp}, $I$ is not
Abelian. Let $J$ be the orthoplement of $I$ with respect to $f$.  As
$f$ is associative, $J$ is an ideal of $L$. We claim  $L = I\oplus J$.
For, $I\cap J \ne0$ would imply $I\subseteq J$ by minimality of $I$
and hence $[I, J]=I$ as $I$ is not Abelian,
so $f(I,L) = f([I,J],L) = f(I,[J,L])\subseteq f(I,J) =0$,
yielding the contradiction $I\subseteq \Rad(f)$.

Hence, if $K$ is a nonzero ideal of $I$, it also is an ideal of $L$,
and so coincides with $I$. This means that $I$ is simple. Moreover, by
Lemma~\ref{orthdecomp}, $J$ is generated by extremal elements with
form $f|_J$.  Notice that $\Rad(f|_J)\subseteq \Rad(f) = 0$, so by
induction on the dimension, we find that $L$ is a direct sum of simple
ideals.
\end{proof}
 
\begin{Cor}
We have that $L/\Rad(L)$ is a direct sum of simple ideals if and only
if $\Rad(L) = \Rad(f)$.
\end{Cor}

\begin{proof}
In Lemma~\ref{ppp}, we have already proved that $\Rad(L) \subseteq \Rad(f)$.
Set $\ov{L} = L/\Rad(L)$ and let $\ov{f}$ be as in
Lemma~\ref{ovf}.

Since $\ov{\Rad(f)} = \Rad(\ov{f})$, we have
$\Rad(L)=\Rad(f)$ if and only if $\Rad(\ov{f}) = 0$, so the corollary
is a direct consequence of Proposition \ref{9.14}.
\end{proof}

\section{Analysis of root groups}
\label{rootgrpsec}

\begin{defn}
For $y \in \Ec$, we define $U_y := \{\exp(y,t) \mid t \in k\}$ to be the
{\rm root group associated to $y$}. Since $\exp(cy,t) =
\exp(y,ct)$ for all $c \in k$, the group only depends on the
1-dimensional subspace $ky$. 
\end{defn}

\bigskip
By $\Sigma := \{ U_y \mid y \in \Ec\}$ we denote the set of root groups
associated to extremal elements of $L$.
Set $G := \erz{\Sigma} \leq \Aut(L)$.
For calculations in $G$, we use  $(\cdot ,\cdot)$ for commutators in their
group theoretic meaning. Thus for $g, h \in G$ and $A, B \leq G$, we
write $(g,h) := g^{-1}h^{-1}gh$ and $(A,B) = \erz{(a,b) \mid a \in A, b \in
B}$. Furthermore, we denote the conjugate of $A$ under $g$ by $A^g $, so $A^g
= g^{-1}Ag$. Observe that $G$ preserves the bilinear form $f$.

\begin{defn}
Let $A$, $B$ be two Abelian subgroups of $G$.
Following Timmesfeld \cite{Tim2}, we call $\erz{A,B}$ a {\em rank 1 group},
if the following holds:
For each $1 \neq a \in A$ there exists some
$1 \neq b \in B$ with $A^b = B^a$; and vice versa.
If in addition $a^b = (b^{-1})^a$, then we call the rank 1 group {\em special}.
\end{defn}

\begin{Thm} \label{abstractroot}
For $x,y \in \Ec$ and $s,t \in k$, the group $G$ has the following
properties.
\begin{enumerate}
\item[(1)] $\exp(y,s)\exp(y,t) = \exp(y,s+t)$. In particular,
$U_y$ is isomorphic to the additive group of $k$. 
\item[(2)] $\exp(y,s) x \in \Ec$ with
$(U_x)^{\exp(y,s)} = U_{\exp(y,-s)x}$.
\item[(3)] If $[x,y] = 0$, then $(U_x,U_y) = 1$.
\item[(4)] If $f(x,y) = 0$, but $[x,y] \neq 0$, then $(\exp(y,t), \exp(x,s)) = \exp([y,x],ts)$.
In particular, the group $\erz{U_x, U_y}$ is nilpotent of class $2$ and
$(U_x,U_y) = (u,U_y) = (U_x,v) = U_{[x,y]}$ for all
$1 \neq u \in U_x$, $1 \neq v \in U_y$. 
\item[(5)] If $f(x,y) = -2$, then, for $s\in k$, $s \neq 0$,
\[\exp(y,-s) \exp(x,s^{-1}t) \exp(y,s) = \exp(x,-s^{-1}) \exp(y,-ts)
\exp(x,s^{-1}).\]
In particular, the group $\erz{U_x, U_y}$ is a special rank
$1$ group in $G$.
\end{enumerate}
\end{Thm}

\begin{proof}
All equations can be verified in a straightforward manner,
by use of the definition of $\exp(y,t)$. In the calculations we use the Jacobi
identity, associativity of $f$
(e.g., $f(y,[y,z]) =  0$, $f(x,[y,[x,z]]) = 0$ when $f(x,y) = 0$,
and $f(x,[y,[x,z]]) = 2f(x,z)$ when $f(x,y) = -2$) and the rewriting rule for
$[x,[y,[x,z]]]$ of Lemma~\ref{ide}. 
\end{proof}

\begin{Remark}
Theorem~\ref{abstractroot} shows that the set $\Sigma$ of root groups
in $G$ associated to extremal elements of $L$ 
is a set of so-called {\em abstract root subgroups}
in the sense of Timmesfeld \cite{Tim2}.
If $\Sigma$ is a conjugacy class and $G$ is quasi-simple, then we may apply
Timmesfeld's classification \cite{Tim2} of groups generated by abstract
root subgroups to determine the structure of $L$.
\end{Remark}

\bigskip
Whenever $G$ is an algebraic group its long root subgroups
correspond to projective points (that is, 1 dimensional linear subspaces)
spanned by extremal elements in its Lie algebra $L(G)$.
The geometry of these long root subgroups, which is well studied,
can thus also be studied in the Lie algebra $L(G)$.

For instance,
in a Chevalley group, consider the subgroup $M$ generated by two different long root
subgroups $U_\alpha$, $ U_\beta$, with $\alpha + \beta$ not a root.
Then either $M$ is partitioned by the long root subgroups contained
in it or it contains no other long root subgroups than $U_\alpha$ and $U_\beta$, 
depending on whether $\alpha - \beta$ is a root or not.
In the lemma below, we describe the corresponding behavior in Lie algebras
for the case of a partitioning. For recognition of such pairs,
we can take the general setup of a Lie algebra containing extremal
elements and do not need to
assume that it is the Lie algebra of an algebraic group.

\begin{Lm} \label{strongcomm}
Let $x, y \in \Ec$ with $[x,y] = 0$. Then the following are equivalent:
\begin{enumerate}
\item[{\rm (1)\phantom{'}}]
For $s,t \in k$, $s,t \neq 0$, the element $sx+ty$ is extremal.
\item[{\rm (1)'}]
There are $s,t \in k$, $s,t \neq 0$, such that
the element $sx+ty$ is extremal.
\item[{\rm (2)\phantom{'}}]
$2[y,[x,z]] = f(x,z)y + f(y,z)x$ for all $z \in \Ec$.
\item[{\rm (2)'}]
$2[y,[x,z]] = f(x,z)y + f(y,z)x$ for all $z \in L$.
\end{enumerate}
Moreover, in these cases, we have
$\exp(y,t)\exp(x,s) = \exp(sx+ty,1)$ for $s,t\in k$.
\end{Lm}

\begin{proof}
Note that (2) and (2)' are equivalent by Lemma~\ref{3.5}.
For $z \in \Ec$, $0 \neq s, t \in k$, we have
$[sx + ty,[sx + ty,z]] 
=  s^2 f(x,z)x + t^2 f(y,z)y + 2st [y,[x,z]]$
and
$f(sx+ty,z)(sx + ty) = s^2f(x,z)x + t^2f(y,z)y +
st\big( f(x,z)y + f(y,z)x \big)$.
Hence (1') implies (2). If (2) holds, then $sx + ty$ is extremal and
a short calculation shows $\exp(y,t)\exp(x,s)z= \exp(sx+ty,1)z$. 
Notice $[x,[y,z]]=[y,[x,z]]$ here.  
This yields the result.
\end{proof}

\begin{Cor}
If three points on a projective line of $L$ represent commuting extremal elements,
then the whole line consists of commuting extremal elements.
\end{Cor}

\bigskip
An example of the occurrence of
projective lines consisting fully of extremal elements
as described in the previous lemma is given in the next lemma.
Note the correspondence with the group geometries of \cite{Tim2}.

\begin{Lm}
Let $x,y \in \Ec$ with $f(x,y) = 0$, but $[x,y] \neq 0$. Then the conditions
of Lemma~{\rm \ref{strongcomm}} hold for $x$ and $[x,y]$.
\end{Lm}

\begin{proof}
By Lemma~\ref{ide} we have $2[[x,y],[x,z]] = f(x,[y,z]) x + f(x,z)[x,y] + 0$.
By associativity of $f$ (cf.\ Theorem \ref{3.6}),
we see $f(x, [y,z]) = f([x,y], z)$, which proves
Condition {\rm (2)'} of Lemma~{\rm \ref{strongcomm}}, as required.
\end{proof}

\bigskip
Many more results in this direction can be derived, such as the non-existence
of a chain $x_1$, $x_2$, $x_3$ of extremal elements such that
$x_1$, $x_2$ are as $x$, $y$
in Lemma \ref{strongcomm},
$[x_2, x_3]=0$,  and $f(x_1,x_3)\ne0$.

To finish, we relate the results on the generation of $L$ 
by extremal elements and to corresponding properties for  $G$.
Let $\mathcal{F}(\Sigma)$ denote the graph whose vertices are the elements
of $\Sigma$ and whose edges are the unordered pairs $\{A,B\}$
with $\erz{A,B}$ a rank $1$ group.

\begin{lemma} \label{genLfromgenG}
Assume that the graph $\mathcal{F}(\Sigma)$
is connected.
If $G = \langle U_{x_i} \mid i \in I \rangle$, then $L$ is generated as a Lie
algebra by the $x_i$ ($i \in I)$. 
\end{lemma}

\begin{proof}
Denote by $W$ the Lie algebra generated by $ \{x_i \mid i \in I\}$.
Let $z \in \Ec$ with $f(z,x_i) \neq 0$ for some $i \in I$.  Since $W$
is invariant under $U_{x_j}$ for $j \in I$, it is also invariant under
$G$.  Hence $[z,x_i] + \frac{1}{2}f(z,x_i)z = \exp(z,1)x_i - x_i \in
W$. We apply $\exp(z,1)$ to this vector and see that $[z,x_i] +
\frac{1}{2}f(z,x_i)z + f(z,x_i)z$ is in $W$. Hence also the
difference, which is $f(z,x_i)z$, is in $W$. Since $f(z,x_i) \neq 0$,
we conclude that $z \in W$.

The graph $\mathcal{F}(\Sigma)$ is connected, so we obtain that
$\Ec \subseteq W$, whence $L = W$.
\end{proof}

\bigskip
For the case where $k$ is algebraically closed and $L$ is of Chevalley type,
we can prove the converse as well.

\begin{Thm}
Suppose that $k$ is algebraically closed.
Let $L$ be a Lie algebra of Chevalley type over $k$ of
characteristic distinct from 2, 3, and let
$G$ be the corresponding group of automorphisms generated by long root groups.
Then $t(L)$, the minimal number of extremal elements generating $L$,
determined in Theorem \ref{mingent}, is the minimal number of
long root groups needed to generate $G$.
\end{Thm}

\begin{proof}
Every element of $\Ec$ is a long root element
by  Proposition \ref{longchev}.
Suppose that $L$
is generated by the extremal elements $x_1,\ldots,x_t$. 
Write  $G=\langle
U_{x} \mid x \in \Ec \rangle$ (as usual) and
 $H = \langle U_{x_i} \mid i=1,\ldots,t\rangle$.
It suffices to show that $H$ coincides with $G$.
Put $\Dc = \{x\in \Ec \mid U_x\subseteq H\}$.
Notice that, for $x,y\in \Dc$, also
$\exp(x,1)y\in \Dc$. Hence, by the argument of the proof of
Lemma \ref{linspan} with $\Dc$ instead of $\Ec$, we obtain that
$L$ is linearly spanned by $\Dc$.

Subgroups of the algebraic group ${\rm GL}(L)$
generated by connected
algebraic subgroups are connected algebraic subgroups
(see \cite{Bor}, Proposition I.2.2),
so $H$ is a closed algebraic subgroup of the connected linear
algebraic group $G$.  Clearly, the
derivative ${\rm d}\iota$ of the embedding $\iota : H\to G$ at the identity
of $H$ is the
embedding $ L(H)\to L(G)$ of the Lie algebra of $H$
in the Lie algebra of $G$.
As $\Dc$ linearly spans $L$, we have $L(H)=L = L(G)$ and so
 ${\rm d}\iota$ is surjective. By
Theorem 3.2.21 of \cite{Spri}, this implies 
that $\iota$ is dominant, that is,
$H$ is dense in $G$. But $H$ is closed as well, so $H=G$.
This establishes that $G$ is generated by at most $t(L)$ long root subgroups.

The converse is handled by the previous lemma.
\end{proof}

\newpage
\bigskip\bigskip
\section*{Addresses of the authors}

\noindent
(correspondence to:)\\
Arjeh M.~Cohen,\\
Fac.~Wisk. en Inf., TUE, POBox 513, \\
5600 MB Eindhoven, The Netherlands\\
(email amc@win.tue.nl)

\medskip\noindent
Anja Steinbach,\\
Math. Inst., Justus-Liebig-Universit\"at,\\
Arndtstr.\ 2, D-35392, Giessen, Germany\\
(email anja.steinbach@math.uni-giessen.de)

\medskip\noindent
Rosane Ushirobira,\\
Fac.~Wisk. en Inf., TUE, POBox 513,\\
5600 MB Eindhoven, The Netherlands\\
(email rosane@win.tue.nl)

\medskip\noindent
David Wales,\\
Sloan Lab,
Caltech,\\
Pasadena, CA 91125, USA\\
(email dbw@caltech.edu)

\end{document}